\documentclass[aos,MSNbibl,nameyear,dvips]{arximspdf}
\usepackage{mathbh}
\usepackage{graphicx}

\doi{10.1214/12-AOS1003} 
\volume{40}
\issue{2}
\pubyear{2012}
\firstpage{1198}
\lastpage{1232}

\makeatletter
\newtheorem{lemma}{Lemma}
\newtheorem{theorem}{Theorem}
\newproclaim{remark}{Remark}

\newcommand{\eqref}[1]{(\ref{#1})}
\newcommand{\argmin}{{\operatorname{argmin}}}
\newcommand{\argmax}{{\operatorname{argmax}}}

\newcommand{\E}{\mathrm{E}}
\renewcommand{\P}{\mathrm{P}}
\newcommand{\Cov}{\operatorname{Cov}}
\newcommand{\half}{\frac{1}{2}}
\newcommand{\tr}{\operatorname{tr}}
\renewcommand{\dim}{\operatorname{dim}}
\newcommand{\col}{\operatorname{col}}
\newcommand{\row}{\operatorname{row}}
\newcommand{\nul}{\operatorname{null}}
\newcommand{\rank}{\operatorname{rank}}
\newcommand{\nuli}{\operatorname{nullity}}
\newcommand{\sign}{\operatorname{sign}}
\newcommand{\supp}{\operatorname{supp}}
\newcommand{\hy}{\hat{y}}
\newcommand{\hbeta}{{\hat{\beta}}}
\newcommand{\hv}{\hat{v}}
\newcommand{\lone}{1}
\newcommand{\ltwo}{2}
\newcommand{\linf}{\infty}
\newcommand{\T}{^T}
\newcommand{\cA}{{\mathcal{A}}}
\newcommand{\cB}{\mathcal{B}}
\newcommand{\cE}{\mathcal{E}}
\newcommand{\cF}{\mathcal{F}}
\newcommand{\cM}{\mathcal{M}}
\newcommand{\cN}{\mathcal{N}}
\newcommand{\cS}{\mathcal{S}}
\newcommand{\df}{\operatorname{df}}
\newcommand{\R}{\mathbb{R}}
\newcommand{\aff}{\operatorname{aff}}
\newcommand{\relint}{\operatorname{relint}}
\newcommand{\relbd}{\operatorname{relbd}}
\makeatother

\begin{document}
\begin{frontmatter}

\title{Degrees of freedom in lasso problems}
\runtitle{Degrees of freedom in lasso problems}

\begin{aug}
\author[A]{\fnms{Ryan J.} \snm{Tibshirani}\corref{}\ead[label=e1]{ryantibs@cmu.edu}}
\and
\author[B]{\fnms{Jonathan} \snm{Taylor}\thanksref{t1}\ead[label=e2]{jtaylo@stanford.edu}}
\runauthor{R. J. Tibshirani and J. Taylor}
\affiliation{Carnegie Mellon University and Stanford University}
\address[A]{Department of Statistics\\
Carnegie Mellon University\\
Baker Hall\\
Pittsburgh, Pennsylvania 15213\\
USA\\
\printead{e1}} 
\address[B]{Department of Statistics\\
            Stanford University\\
            Sequoia Hall\\
            Stanford, California 94305\\
            USA\\
            \printead{e2}}
\end{aug}
\thankstext{t1}{Supported by NSF Grant
DMS-09-06801 and AFOSR Grant 113039.}

\received{\smonth{11} \syear{2011}}
\revised{\smonth{2} \syear{2012}}

%
\begin{abstract}
We derive the degrees of freedom of the lasso fit, placing no
assumptions on the predictor matrix $X$. Like the well-known result of
Zou, Hastie and Tibshirani [\textit{Ann. Statist.}
\textbf{35} (2007) 2173--2192], which gives the degrees of freedom of
the lasso
fit when $X$ has full column rank, we express our result in terms
of the active set of a lasso solution. We extend this result to cover
the degrees of freedom of the generalized lasso fit for an arbitrary
predictor matrix $X$ (and an arbitrary penalty matrix $D$). Though our
focus is degrees of freedom, we establish some intermediate results on
the lasso and generalized lasso that may be interesting on their own.
\end{abstract}

%

\end{frontmatter}

\section{Introduction}
\label{sec:intro}

We study degrees of freedom, or the ``effective number of
parameters,'' in $\ell_1$-penalized linear regression problems. In
particular, for a response vector $y \in\R^n$, predictor matrix $X
\in\R^{n \times p}$ and tuning parameter $\lambda\geq0$, we
consider the lasso problem [\citet{bp}, \citet{lasso}]
%
\begin{equation}
\label{eq:lasso}
\hbeta\in\mathop\argmin_{\beta\in\R^p} \half\|y-X\beta\|_\ltwo^2 +
\lambda\|\beta\|_\lone.
\end{equation}
The above notation emphasizes the fact that the solution $\hbeta$ may
not be unique [such nonuniqueness can occur if
$\rank(X)<p$]. Throughout the paper, when a function $f\dvtx D
\rightarrow\R^n$ may have a nonunique minimizer over its domain $D$,
we write $\argmin_{x \in D} f(x)$ to denote the set of minimizing $x$
values, that is, $\argmin_{x \in D} f(x) = \{ \hat{x}
\in D \dvtx f(\hat{x}) = \min_{x \in D} f(x) \}$.

A fundamental result on the degrees of freedom of the lasso fit was
shown by \citet{lassodf}. The authors show that if $y$ follows
a~normal distribution with spherical\vadjust{\goodbreak} covariance,
$y \sim N(\mu,\sigma^2 I)$, and $X,\lambda$ are considered fixed with
$\rank(X)=p$, then
%
\begin{equation}
\label{eq:lassodffull}
\df(X\hbeta) = \E|\cA|,
\end{equation}
where $\cA=\cA(y)$ denotes the active set of the unique lasso
solution at $y$, and~$|\cA|$ is its cardinality. This is quite a
well-known result, and is sometimes used to informally justify an
application of the lasso procedure, as it says that
number of parameters used by the lasso fit is simply equal to the
(average) number of selected variables. However, we note that
the assumption $\rank(X)=p$ implies that $p \leq n$; in other words,
the degrees of freedom result \eqref{eq:lassodffull} does not cover
the important ``high-dimensional'' case $p>n$. In this case,
the lasso solution is not necessarily unique, which raises
the questions:
\begin{itemize}
\item Can we still express degrees of freedom in terms of the
active set of a lasso solution?
\item If so, which active set (solution) would we refer to?
\end{itemize}
In Section~\ref{sec:lasso}, we provide answers to these questions, by
proving a stronger result when $X$ is a general predictor matrix.
We show that the subspace spanned by the columns of $X$ in
$\cA$ is almost surely unique, where ``almost surely'' means for
almost every $y \in\R^n$. Furthermore, the degrees of freedom of the
lasso fit is simply the expected dimension of this column space.

We also consider the generalized lasso problem,
%
\begin{equation}
\label{eq:genlasso}
\hbeta\in\mathop\argmin_{\beta\in\R^p} \half\|y-X\beta\|_\ltwo^2 +
\lambda\|D\beta\|_\lone,
\end{equation}
where $D \in\R^{m \times p}$ is a penalty matrix, and again the
notation emphasizes the fact that $\hbeta$ need not be unique
[when $\rank(X)<p$]. This of course reduces to the usual
lasso problem \eqref{eq:lasso} when $D=I$, and \citet{genlasso}
demonstrate that the formulation \eqref{eq:genlasso} encapsulates
several other important problems---including the fused lasso on any
graph and trend filtering of any order---by varying the penalty matrix
$D$. The same paper shows that if
$y$ is normally distributed as above, and $X,D,\lambda$ are fixed with
$\rank(X)=p$, then the generalized lasso fit has degrees of freedom
%
\begin{equation}
\label{eq:genlassodffull}
\df(X\hbeta) = \E[\nuli(D_{-\cB})].
\end{equation}
Here $\cB=\cB(y)$ denotes the boundary set of an optimal
subgradient to the generalized lasso problem at $y$ (equivalently, the
boundary set of a dual solution at $y$), $D_{-\cB}$ denotes the matrix $D$
after having removed the rows that are indexed by $\cB$, and
$\nuli(D_{-\cB})=\dim(\nul(D_{-\cB}))$, the dimension of the null
space of
$D_{-\cB}$.

It turns out that examining \eqref{eq:genlassodffull} for specific
choices of $D$ produces a number of interpretable corollaries, as
discussed in \citet{genlasso}. For
example, this result implies that the degrees of freedom of the fused
lasso fit is equal to the expected number of fused groups, and that
the degrees of freedom of the\vadjust{\goodbreak}
trend filtering fit is equal to the expected number of knots
$+$ $(k-1)$, where $k$ is the order of the polynomial.
The result \eqref{eq:genlassodffull} assumes that $\rank(X)=p$ and
does not cover the case $p>n$; in Section~\ref{sec:genlasso},
we derive the degrees of freedom of the generalized lasso fit for a
general $X$ (and still a~general $D$).
As in the lasso case, we prove that there exists a linear subspace
$X(\nul(D_{-\cB}))$ that is almost surely unique, meaning that it will
be the same under different boundary sets $\cB$ corresponding to
different solutions of \eqref{eq:genlasso}. The generalized lasso
degrees of freedom is then the expected dimension of this subspace.

Our assumptions throughout the paper are minimal. As was already
mentioned, we place no assumptions whatsoever on the predictor matrix
$X \in\R^{n \times p}$ or on the penalty matrix $D \in\R^{m \times p}$,
considering them fixed and nonrandom. We also consider $\lambda\geq
0$ fixed. For Theorems~\ref{thm:lassodfequi},~\ref{thm:lassodfact}
and~\ref{thm:genlassodf} we assume that $y$ is normally distributed,
%
\begin{equation}
\label{eq:normal}
y \sim N(\mu,\sigma^2 I)
\end{equation}
for some (unknown) mean vector $\mu\in\R^n$ and marginal variance
$\sigma^2 \geq0$. This assumption is only needed in order to apply
Stein's formula for degrees of freedom, and none of the other lasso
and generalized lasso results in the paper, namely
Lemmas~\ref{lem:lassoproj} through~\ref{lem:invbound}, make any assumption
about the distribution of $y$.

This paper is organized as follows. The rest of the Introduction
contains an overview of related work, and an explanation of our notation.
Section~\ref{sec:prelim} covers some relevant background material on
degrees of freedom and convex polyhedra. Though the connection may not
be immediately obvious, the geometry of polyhedra plays a large
role in understanding problems \eqref{eq:lasso} and
\eqref{eq:genlasso}, and Section~\ref{sec:poly} gives a high-level
view of this geometry before the technical arguments
that follow in Sections~\ref{sec:lasso} and~\ref{sec:genlasso}.
In Section~\ref{sec:lasso}, we derive two representations
for the degrees of freedom of the lasso fit, given in
Theorems~\ref{thm:lassodfequi} and~\ref{thm:lassodfact}. In
Section~\ref{sec:genlasso}, we derive the analogous results for the
generalized lasso problem, and these are given in
Theorem~\ref{thm:genlassodf}. As the lasso problem is a special case of the
generalized lasso problem (corresponding to $D=I$),
Theorems~\ref{thm:lassodfequi} and~\ref{thm:lassodfact} can actually be viewed
as corollaries of Theorem~\ref{thm:genlassodf}. The reader may then
ask: why is there a separate section dedicated to the lasso
problem? We give two reasons: first, the lasso arguments are simpler
and easier to follow than their generalized lasso counterparts;
second, we cover some intermediate results for the lasso problem that
are interesting in their own right and that do not carry over to the
generalized lasso perspective. Section~\ref{sec:disc} contains some
final discussion.

\subsection{Related work}
\label{sec:work}
All of the degrees of freedom results discussed here assume that
the response vector has distribution $y \sim N(\mu,\sigma^2 I)$, and that
the predictor matrix $X$ is fixed.
To the best of our knowledge, \citet{lars} were the first to prove a
result on
the degrees of freedom of the lasso fit, using the lasso solution path
with $\lambda$ moving from $\infty$ to $0$.
The authors showed that when the active set
reaches size $k$ along this path, the lasso fit has degrees of freedom exactly
$k$. This result assumes that $X$ has full column rank and further
satisfies a
restrictive condition called the ``positive cone condition,'' which
ensures that
as $\lambda$ decreases, variables can only enter, and not leave, the
active set.
Subsequent results on the lasso degrees of freedom (including those
presented in this
paper) differ from this original result in that they derive degrees of
freedom for
a fixed value of the tuning parameter $\lambda$, and not a fixed
number of steps
$k$ taken along the solution path.

As mentioned previously, \citet{lassodf} established the basic lasso
degrees of freedom
result (for fixed $\lambda$) stated in \eqref{eq:lassodffull}. This
is analogous
to the path result of \citet{lars}; here degrees of freedom is equal to the
expected size of the active set (rather than simply the size) because
for a fixed
$\lambda$ the active set is a random quantity, and can hence achieve a
random size.
The proof of \eqref{eq:lassodffull} appearing in \citet{lassodf}
relies heavily on
properties of the lasso solution path. As also mentioned previously,
\citet{genlasso}
derived an extension of~\eqref{eq:lassodffull} to the generalized
lasso problem, which is stated in \eqref{eq:genlassodffull} for an
arbitrary penalty matrix
$D$. Their arguments are not based on properties of the solution path,
but instead
come from a geometric perspective much like the one developed in this paper.

Both of the results \eqref{eq:lassodffull} and \eqref
{eq:genlassodffull} assume that
$\rank(X)=p$; the current work extends these to the case of an
arbitrary matrix $X$, in
Theorems~\ref{thm:lassodfequi},~\ref{thm:lassodfact} (the lasso)
and~\ref{thm:genlassodf}
(the generalized lasso). In terms of our intermediate results, a
version of Lemmas~\ref{lem:lcequi},~\ref{lem:lcact} corresponding to $\rank(X)=p$
appears in \citet{lassodf},
and a version of Lemma~\ref{lem:lcbound} corresponding to $\rank
(X)=p$ appears in
\citet{genlasso} [furthermore, \citet{genlasso} only consider the
boundary set
representation and not the active set representation]. Lemmas~\ref{lem:nonexp},
\ref{lem:locaff} and the conclusions thereafter, on the degrees of
freedom of the projection
map onto a convex polyhedron, are essentially given in \citet
{meyerwood}, though these authors
state and prove the results in a different manner.

In preparing a draft of this manuscript, it was brought to our
attention that other
authors have independently and concurrently worked to extend results
\eqref{eq:lassodffull}
and \eqref{eq:genlassodffull} to the general $X$ case. Namely, \citet
{dfl1} prove
a result on the lasso degrees of freedom, and \citet{analreg} prove a~result on the
generalized lasso degrees of freedom, both for an arbitrary $X$. These
authors' results
express degrees of freedom in terms of the active sets of special
(lasso or generalized
lasso) solutions. Theorems~\ref{thm:lassodfact} and~\ref
{thm:genlassodf} express degrees
of freedom in terms of the active sets of any solutions, and hence the
appropriate application
of these theorems provides an alternative verification of these
formulas. We discuss this in
detail in the form of remarks following the theorems.

\subsection{Notation}
\label{sec:notat}

In this paper, we use $\col(A)$, $\row(A)$ and $\nul(A)$ to denote the
column space, row space and null space of a matrix $A$, respectively;
we use $\rank(A)$ and $\nuli(A)$ to denote the dimensions of $\col(A)$
[equivalently, $\row(A)$] and $\nul(A)$, respectively. We write $A^+$
for the the Moore--Penrose pseudoinverse of $A$; for a rectangular
matrix $A$, recall that $A^+=(A\T A)^+ A\T$. We write $P_L$ to denote
the projection matrix onto a linear subspace $L$, and more
generally, $P_C(x)$ to denote the projection of a~point~$x$ onto a
closed convex set $C$. For readability, we sometimes write
$\langle a, b \rangle$ (instead of $a\T b$) to denote the inner
product between vectors $a$ and $b$.

For a set of indices $R = \{i_1,\ldots, i_k\} \subseteq\{1,\ldots,
m\}$ satisfying $i_1 < \cdots< i_k$, and a vector $x \in\R^m$, we
use $x_R$ to denote the subvector $x_R = (x_{i_1},\ldots, x_{i_k})\T
\in\R^k$. We denote the complementary subvector by
$x_{-R}=x_{\{1, \ldots, m\}\setminus R} \in\R^{m-k}$.
The notation is similar for matrices.
Given another subset of indices $S = \{j_1,\ldots, j_\ell\}
\subseteq\{1,\ldots, p\}$ with $j_1
< \cdots< j_\ell$, and a matrix $A \in\R^{m \times p}$, we use
$A_{(R,S)}$ to denote the submatrix
\[
A_{(R,S)} = \left[\matrix{
A_{i_1,j_1} & \cdots& A_{i_1,j_\ell}
\cr
\vdots& &
\cr
A_{i_k,j_1} & \cdots& A_{i_k,j_\ell}}
\right] \in\R^{k \times\ell}.
\]
In words, rows are indexed by $R$, and columns are indexed by
$S$. When combining this notation with the transpose operation, we
assume that the indexing happens first, so that
$A_{(R,S)}\T=(A_{(R,S)})\T$. As above, negative signs are used
to denote the complementary set of rows or columns; for example,
$A_{(-R,S)}=A_{(\{1, \ldots, m\}\setminus R, S)}$. To extract only
rows or only columns, we abbreviate the other dimension by a dot,
so that $A_{(R,\cdot)}=A_{(R,\{1,\ldots, p\})}$ and
$A_{(\cdot,S)}=A_{(\{1,\ldots, m\},S)}$; to extract a single
row or column, we use $A_{(i,\cdot)}=A_{(\{i\},\cdot)}$ or
$A_{(\cdot,j)}=A_{(\cdot,\{j\})}$.
Finally, and most importantly, we introduce the following shorthand
notation:
\begin{itemize}
\item For the predictor matrix $X \in\R^{n\times p}$, we let
$X_S=X_{(\cdot,S)}$.
\item For the penalty matrix $D \in\R^{m \times p}$, we let
$D_R=D_{(R,\cdot)}$.
\end{itemize}
In other words, the default for $X$ is to index its columns, and the
default for $D$ is to index its rows. This convention greatly
simplifies the notation in expressions that involve multiple instances
of $X_S$ or $D_R$; however, its use could also cause a great deal of
confusion, if not properly interpreted by the reader!

\section{Preliminary material}
\label{sec:prelim}

The following two sections describe some background material needed
to follow the results in Sections~\ref{sec:lasso} and~\ref{sec:genlasso}.

\subsection{Degrees of freedom}
\label{sec:df}

If the data vector $y \in\R^n$ is distributed
according to the homoskedastic model
$y \sim(\mu,\sigma^2 I)$,
meaning that the components of $y$ are uncorrelated, with $y_i$
having mean $\mu_i$ and variance $\sigma^2$ for $i=1,\ldots,
n$, then the degrees of freedom of a function\vadjust{\goodbreak}
$g\dvtx \R^n \rightarrow\R^n$ with
$g(y)=(g_1(y),\ldots, g_n(y))\T$, is defined as
%
\begin{equation}
\label{eq:df}
\df(g) = \frac{1}{\sigma^2} \sum_{i=1}^n \Cov
(g_i(y),y_i ).
\end{equation}
This definition is often attributed to \citet{bradbiased} or
\citet{gam}, and is interpreted as the ``effective number of
parameters'' used by the fitting procedure $g$. Note that for the
linear regression fit of $y \in\R^n$ onto a fixed and full
column rank predictor matrix $X \in\R^{n \times p}$, we have $g(y) =
\hy= XX^+y$, and $\df(\hy)=\tr(XX^+) = p$, which is the number of
fitted coefficients (one for each predictor variable). Furthermore, we
can decompose the risk of $\hy$, denoted by $\operatorname{Risk}(\hy) =
\E\|\hy-\mu\|_\ltwo^2$, as
\[
\operatorname{Risk}(\hy) = \E\|\hy-y\|_\ltwo^2-
n\sigma^2 + 2p\sigma^2,
\]
a well-known identity that leads to the derivation of the $C_p$
statistic [\citet{mallows}]. For a general fitting procedure $g$, the
motivation for the definition \eqref{eq:df} comes from the analogous
decomposition of the quantity $\operatorname{Risk}(g) =
\E\|g(y)-\mu\|_\ltwo^2$,
%
\begin{equation}
\label{eq:riskd}
\operatorname{Risk}(g) = \E\|g(y)-y\|_\ltwo^2 -
n\sigma^2 + 2\sum_{i=1}^n\Cov(g_i(y),y_i ).
\end{equation}
Therefore a large difference between risk and expected training error
implies a large degrees of freedom.

Why is the concept of degrees of freedom important? One simple answer
is that it provides a way to put different fitting procedures on
equal footing. For example, it would not seem fair to compare a procedure
that uses an effective number of parameters equal to 100 with another that
uses only 10. However, assuming that these procedures can be tuned
to varying levels of adaptivity (as is the case with the lasso and
generalized lasso, where the
adaptivity is controlled by $\lambda$), one could first tune the procedures
to have the same degrees of freedom, and then compare their
performances. Doing this over several common values for degrees of freedom
may reveal, in an informal sense, that one procedure is particularly
efficient when it comes to its parameter usage versus another.

A more detailed answer to the above question is based the risk
decomposition~\eqref{eq:riskd}. The decomposition suggests that an
estimate $\widehat{\df}(g)$ of degrees of freedom can be used
to form an estimate of the risk,
%
\begin{equation}
\label{eq:riskhat}
\widehat{\operatorname{Risk}}(g) = \|g(y)-y\|_\ltwo^2 -n\sigma^2 +
2\sigma^2 \widehat{\df}(g).
\end{equation}
Furthermore, it is straightforward to check that an unbiased
estimate of degrees of freedom leads to an unbiased estimate of risk;
that is, $\df(g)=\E[\widehat{\df}(g)]$ implies $\operatorname{Risk}(g) =
\E[\widehat{\operatorname{Risk}}(g)]$. Hence, the risk estimate \eqref{eq:riskhat}
can be used to choose between fitting procedures, assuming that unbiased
estimates of degrees of freedom are available.
[It is worth mentioning that bootstrap or Monte Carlo methods can\vadjust{\goodbreak}
be helpful in estimating degrees of freedom~\eqref{eq:df} when an analytic
form is difficult to obtain.]
The natural extension of this idea is to use the risk estimate
\eqref{eq:riskhat} for tuning parameter selection. If we suppose that $g$
depends on a tuning parameter $\lambda\in\Lambda$, denoted
$g=g_\lambda(y)$,
then in principle one could minimize the estimated risk over $\lambda$ to
select an appropriate value for the tuning parameter,
%
\begin{equation}
\label{eq:tunsel}
\hat{\lambda} = \mathop\argmin_{\lambda\in\Lambda} \|g_\lambda(y)-y\|
_\ltwo^2 -
n\sigma^2 + 2\sigma^2 \widehat{\df}(g_\lambda).
\end{equation}
This is a computationally efficient alternative to selecting the tuning
parameter by cross-validation, and it is commonly used (along with similar
methods that replace the factor of $2$ above
with a function of $n$ or $p$) in penalized regression problems.
Even though such an estimate \eqref{eq:tunsel} is commonly used in the
high-dimensional setting ($p>n$), its asymptotic properties are
largely unknown in this case, such as risk consistency, or relatively
efficiency compared to the cross-validation estimate.

\citet{stein} proposed the risk estimate \eqref{eq:riskhat} using a
particular unbiased estimate of degrees of freedom, now
commonly referred to as \textit{Stein's unbiased risk estimate} (SURE).
Stein's framework requires that we strengthen our distributional
assumption on $y$ and assume normality, as stated in
\eqref{eq:normal}. We also assume that the function $g$ is continuous
and almost differentiable. (The precise definition of almost
differentiability is not important here, but the interested reader may
take it to mean that each coordinate function~$g_i$ is absolutely
continuous on almost every line segment parallel to one of the coordinate
axes.) Given these assumptions, Stein's main result is an alternate
expression for degrees of freedom,
%
\begin{equation}
\label{eq:steindf}
\df(g) = \E[(\nabla\cdot g)(y)],
\end{equation}
where the function $\nabla\cdot g =
\sum_{i=1}^n \partial g_i/\partial y_i$
is called the divergence of $g$. Immediately following is the
unbiased estimate of degrees of freedom,
%
\begin{equation}
\label{eq:steindfhat}
\widehat{\df}(g) = (\nabla\cdot g)(y).
\end{equation}
We pause for a moment to reflect on the importance of this result.
From its definition \eqref{eq:df}, we can see that the two most
obvious candidates for unbiased estimates of degrees of freedom are
\[
\frac{1}{\sigma^2} \sum_{i=1}^n g_i(y)(y_i-\mu_i)
\quad\mbox{and}\quad
\frac{1}{\sigma^2} \sum_{i=1}^n \bigl(g_i(y)-\E[g_i(y)] \bigr) y_i.
\]
To use the first estimate above, we need to know $\mu$ (remember, this
is ultimately what we are trying to estimate!). Using
the second requires knowing $\E[g(y)]$, which is equally impractical
because this invariably depends on $\mu$. On the other hand, Stein's
unbiased estimate \eqref{eq:steindfhat} does not have an explicit
dependence on $\mu$; moreover, it can be analytically computed for
many fitting procedures $g$. For example, Theorem~\ref{thm:lassodfact} in
Section~\ref{sec:lasso} shows that, except for $y$ in a set of measure
zero, the divergence of the lasso fit is equal to $\rank(X_\cA)$ with\vadjust{\goodbreak}
$\cA=\cA(y)$ being the active set of a lasso solution at $y$.
Hence, Stein's formula allows for the unbiased estimation of
degrees of freedom (and subsequently, risk) for a broad class of
fitting procedures $g$---something that may have not seemed possible
when working from the definition directly.

\subsection{Projections onto polyhedra}
\label{sec:poly}

\begin{figure}

\includegraphics{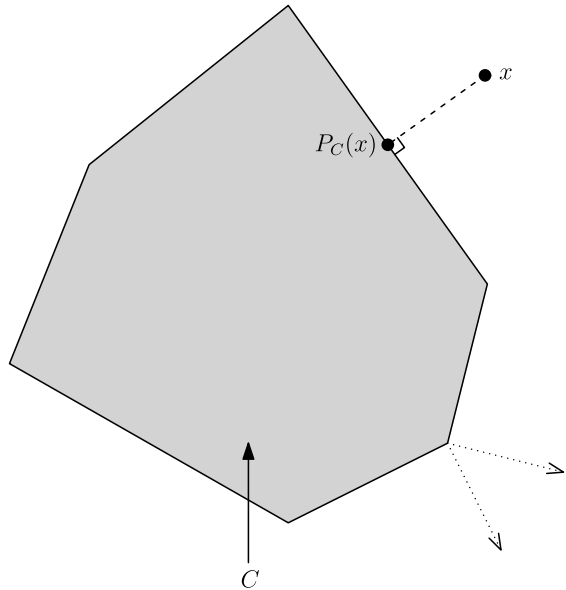}

\caption{\it An example of a polyhedron in
$\R^2$.}\label{fig:poly}
\end{figure}

A set $C \subseteq\R^n$ is called a \textit{convex polyhedron}, or
simply a \textit{polyhedron}, if $C$ is the intersection of finitely
many half-spaces,
%
\begin{equation}
\label{eq:poly}
C = \bigcap_{i=1}^k \{x \in\R^n \dvtx a_i\T x \leq b_i \},
\end{equation}
where $a_1,\ldots, a_k \in\R^n$ and $b_1,\ldots, b_k \in\R$.
(Note that we do not require boundedness here; a
bounded polyhedron is sometimes called a polytope.) See
Figure~\ref{fig:poly} for an example. There is a rich theory on
polyhedra; the definitive reference is \citet{grunbaum}, and
another good reference is \citet{schneider}. As this
is a paper on statistics and not geometry, we do not attempt to give
an extensive treatment of the properties of polyhedra. We do,
however, give two properties (in the form of two lemmas) that are
especially important with respect to our statistical problem; our
discussion will also make it clear why polyhedra are relevant
in the first place.

From its definition \eqref{eq:poly}, it follows that a polyhedron
is a closed convex set. The first property that we discuss does not
actually rely on the special structure of polyhedra, but only on
convexity. For any closed convex set
$C \subseteq\R^n$ and any point $x \in\R^n$, there is a unique point
$u \in C$ minimizing $\|x-u\|_\ltwo$. To see this, note that if $v \in
C$ is another minimizer, $v \not= u$, then by convexity $w = (u+v)/2
\in C$, and $\|x-w\|_\ltwo< \|x-u\|_\ltwo/2 + \|x-v\|_\ltwo/2 =
\|x-u\|_\ltwo$, a~contradiction. Therefore, the projection map onto
$C$ is indeed well defined, and we write this as $P_C \dvtx \R^n
\rightarrow C$,
\[
P_C(x) = \mathop\argmin_{u \in C} \|x-u\|_\ltwo.
\]
For the usual linear regression problem, where $y \in\R^n$ is
regressed onto $X \in\R^{n \times p}$, the fit $X\hbeta$ can be
written in terms of the projection map onto the polyhedron
$C=\col(X)$, as in $X\hbeta(y) = X X^+ y =
P_{\col(X)}(y)$. Furthermore, for both the lasso
and generalized lasso problems, \eqref{eq:lasso} and
\eqref{eq:genlasso}, it turns out that we can express the fit as the
residual from projecting onto a suitable polyhedron $C \subseteq
\R^n$, that is,
\[
X\hbeta(y) = (I-P_C)(y) = y-P_C(y).
\]
This is proved in Lemma~\ref{lem:lassoproj} for the lasso and in
Lemma~\ref{lem:genlassoproj} for the generalized lasso (the polyhedron
$C$ depends on $X,\lambda$ for the lasso case, and on $X,D,\lambda$
for the generalized lasso case). Our first lemma establishes that both
the projection map onto a closed convex set and the residual map are
nonexpansive, hence continuous and almost differentiable
everywhere. These are the conditions needed to apply Stein's formula.

\begin{lemma}
\label{lem:nonexp}
For any closed convex set $C \subseteq\R^n$, both the projection map
$P_C \dvtx \R^n \rightarrow C$ and the residual projection map $I-P_C
\dvtx
\R^n \rightarrow\R^n$ are nonexpansive. That is, they satisfy
\begin{eqnarray*}
\|P_C(x) - \P_C(y)\|_\ltwo
&\leq&\|x-y\|_\ltwo\qquad
\mbox{for any } x,y \in\R^n, \quad\mbox{and}
\\
\|(I-P_C)(x) - (I-P_C)(y)\|_\ltwo
&\leq&\|x-y\|_\ltwo\qquad
\mbox{for any } x,y \in\R^n.
\end{eqnarray*}
Therefore, $P_C$ and $I-P_C$ are both continuous and almost
differentiable.
\end{lemma}

The proof can be found in Appendix~\ref{app:nonexp}.
Lemma~\ref{lem:nonexp} will be quite useful later in the paper, as it
will allow us to use Stein's formula to compute the degrees of freedom
of the lasso and generalized lasso fits, after showing that these
fits are indeed the residuals from projecting onto closed convex sets.

The second property that we discuss uses the structure of
polyhedra. Unlike Lemma~\ref{lem:nonexp}, this property will not be
used directly in the following sections of the paper; instead, we
present it here to give some intuition with respect to the degrees of
freedom calculations to come. The property can be best explained by
looking back at Figure~\ref{fig:poly}. Loosely speaking, the picture
suggests that we can move the point $x$ around a bit and it will still
project to the same face of $C$. Another way of saying this is that
there is a neighborhood of $x$ on which $P_C$ is simply the
projection onto an affine subspace. This would not be true if $x$ is
in some exceptional set, which is made up of rays that emanate from
the corners of $C$, like the two drawn in the bottom right corner
of figure. However, the union of such rays has measure zero, so the
map $P_C$ is locally an affine projection, almost everywhere. This
idea can be stated formally as follows.

\begin{lemma}
\label{lem:locaff}
Let $C \subseteq\R^n$ be a polyhedron. For almost every $x \in
\R^n$, there is an associated neighborhood $U$ of $x$, linear
subspace $L \subseteq\R^n$ and point $a \in\R^n$, such that the
projection map restricted to $U$, $P_C\dvtx U \rightarrow C$, is
\[
P_C(y) = P_L(y-a)+a \qquad\mbox{for } y \in U,
\]
which is simply the projection onto the affine subspace $L+a$.
\end{lemma}

The proof is given in Appendix~\ref{app:locaff}.
These last two properties can be used to derive a
general expression for the degrees of freedom of the fitting procedure
$g(y)=(I-P_C)(y)$, when $C \subseteq\R^n$ is a polyhedron.
[A similar formula holds for $g(y)=P_C(y)$.]
Lemma~\ref{lem:nonexp} tells us that $I-P_C$ is continuous and almost
differentiable, so we can use Stein's formula \eqref{eq:steindf}
to compute its degrees of
freedom. Lemma~\ref{lem:locaff} tells us that for almost every $y \in
\R^n$, there is a neighborhood $U$ of $y$, linear subspace $L
\subseteq\R^n$, and point $a \in\R^n$, such that
\[
(I-P_C)(y') = y'-P_L(y'-a)-a = (I-P_L)(y'-a)
\qquad\mbox{for } y' \in U.
\]
Therefore,
\[
\bigl(\nabla\cdot(I-P_C)\bigr)(y) = \tr(I-P_L) = n-\dim(L),
\]
and an expectation over $y$ gives
\[
\df(I-P_C)=n-\E[\dim(L)].
\]
It should be made clear that the random quantity in the above
expectation is the linear subspace $L=L(y)$, which depends on $y$.

In a sense, the remainder of this paper is focused on describing
$\dim(L)$---the dimension of the face of $C$ onto which the point $y$
projects---in a meaningful way for the lasso and generalized lasso
problems. Section~\ref{sec:lasso} considers the lasso problem, and we
show that $L$ can be written in terms of the equicorrelation set of
the fit at $y$. We also show that $L$ can be described
in terms of the active set of a solution at $y$. In
Section~\ref{sec:genlasso} we show the analogous results for the generalized
lasso problem, namely, that $L$ can be written in terms of either the
boundary set of an optimal subgradient at $y$ (the analogy of the
equicorrelation set for the lasso) or the active set of a solution at
$y$.

\section{The lasso}
\label{sec:lasso}

In this section we derive the degrees of freedom of the lasso fit, for
a general predictor matrix $X$. All of our arguments stem from the
Karush--Kuhn--Tucker (KKT) optimality conditions, and we present these
first.\vadjust{\goodbreak}
We note that many of the results in this section can be alternatively
derived using the lasso dual problem.
Appendix~\ref{app:dual} explains this connection more precisely. For the current
work, we avoid the dual perspective simply to keep the
presentation more self-contained. Finally, we
remind the reader that $X_S$ is used to extract columns of $X$
corresponding to an index set~$S$.\looseness=1

\subsection{The KKT conditions and the underlying polyhedron}
\label{sec:lassokkt}

The KKT conditions for the lasso problem \eqref{eq:lasso} can be
expressed as
%
\begin{eqnarray}
\label{eq:lassokkt}
&\displaystyle X\T(y-X\hbeta) = \lambda\gamma,&
\\
\label{eq:lassosg}
&\displaystyle\gamma_i \in
\cases{
\{\sign(\hbeta_i)\} &\quad if  $\hbeta_i \not= 0$,\vspace*{2pt}
\cr
[-1,1] & \quad if $\hbeta_i = 0$.}&
\end{eqnarray}
Here $\gamma\in\R^p$ is a subgradient of the function $f(x) =
\|x\|_\lone$ evaluated at $x=\hbeta$. Hence $\hbeta$ is a
minimizer in \eqref{eq:lasso} if and only if $\hbeta$ satisfies
\eqref{eq:lassokkt} and \eqref{eq:lassosg} for some~$\gamma$.
Directly from the KKT conditions, we can show that $X\hbeta$ is
the residual from projecting $y$ onto a polyhedron.

\begin{lemma}
\label{lem:lassoproj}
For any $X$ and $\lambda\geq0$, the lasso fit
$X\hbeta$ can be written as $X\hbeta(y) = (I-P_C)(y)$, where $C
\subseteq\R^n$ is the polyhedron
\[
C = \{ u \in\R^n \dvtx \|X\T u\|_\linf\leq\lambda\}.
\]
\end{lemma}
\begin{pf}
Given a point $y \in\R^n$, its projection $\theta= P_C(y)$
onto a closed convex set $C \subseteq\R^n$ can be characterized as
the unique point satisfying
%
\begin{equation}
\label{eq:opt}
\langle y-\theta, \theta-u \rangle\geq0\qquad \mbox{for all }
u \in C.
\end{equation}
Hence defining $\theta=y-X\hbeta(y)$, and $C$ as in the lemma,
we want to show that~\eqref{eq:opt} holds for all $u \in C$.
Well,
%
\begin{eqnarray}\label{eq:innerprod}
\langle y-\theta, \theta-u\rangle
&=& \langle X\hbeta, y-X\hbeta- u\rangle\nonumber
\\[-8pt]
\\[-8pt]
&=& \langle X\hbeta, y-X\hbeta\rangle
- \langle X\T u, \hbeta\rangle.\nonumber
\end{eqnarray}
Consider the first term above. Taking an inner product with $\hbeta$
on both sides of \eqref{eq:lassokkt} gives $\langle X\hbeta,
y-X\hbeta\rangle= \lambda\|\hbeta\|_\lone$. Furthermore, the
$\ell_1$ norm can be characterized in terms of its dual norm, the
$\ell_\infty$ norm, as in
\[
\lambda\|\hbeta\|_\lone= \max_{\|w\|_\linf\leq\lambda}
\langle w, \hbeta\rangle.
\]
Therefore, continuing from \eqref{eq:innerprod}, we have
\[
\langle y-\theta, \theta-u \rangle=
\max_{\|w\|_\linf\leq\lambda} \langle w, \hbeta\rangle
- \langle X\T u, \hbeta\rangle,
\]
which is $\geq0$ for all $u \in C$, and we have hence proved
that $\theta=y-X\hbeta(y) = P_C(y)$. To show that $C$ is indeed a
polyhedron, note that it can be written~as
\[
C = \bigcap_{i=1}^p (
\{u \in\R^n \dvtx X_i\T u \leq\lambda\} \cap
\{u \in\R^n \dvtx X_i\T u \geq-\lambda\} ),
\]
which is a finite intersection of half-spaces.
\end{pf}

Showing that the lasso fit is the residual from projecting $y$ onto a
polyhedron is important, because it means that
$X\hbeta(y)$ is nonexpansive as a~function of $y$, and hence
continuous and almost differentiable, by Lemma~\ref{lem:nonexp}. This
establishes the conditions that are needed to apply Stein's formula
for degrees of freedom.

In the next section, we define the equicorrelation set $\cE$, and show
that the lasso fit and solutions both have an explicit form in terms
of $\cE$. Following this, we derive an expression for the
lasso degrees of freedom as a function of the equicorrelation set.

\subsection{The equicorrelation set}
\label{sec:lassoequi}

According to Lemma~\ref{lem:lassoproj}, the lasso fit $X\hbeta$ is
always unique (because projection onto a closed convex set is
unique).\vspace*{1pt} Therefore, even though the solution $\hbeta$ is not
necessarily unique, the optimal subgradient $\gamma$ is unique,
because it can be written entirely in terms of $X\hbeta$, as shown
by \eqref{eq:lassokkt}. We define the unique \textit{equicorrelation set}
$\cE$ as
%
\begin{equation}
\label{eq:equiset}
\cE= \{ i \in\{1,\ldots, p\} \dvtx |\gamma_i| = 1 \}.
\end{equation}
An alternative definition for the equicorrelation set is
%
\begin{equation}
\label{eq:equiset2}
\cE= \{i \in\{1,\ldots, p\} \dvtx |X_i \T(y-X\hbeta)| =
\lambda\},
\end{equation}
which explains its name, as $\cE$ can be thought of as the set of
variables that have equal and maximal absolute inner product
(or correlation for standardized variables) with the residual.

The set $\cE$ is a natural quantity to work with because we can
express the lasso fit and the set of lasso solutions in terms of
$\cE$, by working directly from equation~\eqref{eq:lassokkt}.
First we let
%
\begin{equation}
\label{eq:equisigns}
s = \sign(\gamma_\cE) = \sign\bigl(X_\cE\T(y-X\hbeta)\bigr),
\end{equation}
the signs of the inner products of the equicorrelation variables with
the residual. Since $\hbeta_{-\cE}=0$ by definition of
the subgradient, the $\cE$ block of the KKT conditions
can be rewritten as
%
\begin{equation}
\label{eq:lassokkt2}
X_\cE\T(y-X_\cE\hbeta_\cE) = \lambda s.
\end{equation}
Because $\lambda s \in\row(X_\cE)$, we can write $\lambda s =
X_\cE\T(X_\cE\T)^+ \lambda s$, so rearranging
\eqref{eq:lassokkt2} we get
\[
X_\cE\T X_\cE\hbeta_\cE=
X_\cE\T \bigl(y-(X_\cE\T)^+ \lambda s \bigr).
\]
Therefore, the lasso fit $X\hbeta= X_\cE\hbeta_\cE$ is
%
\begin{equation}
\label{eq:lassofit}
X\hbeta= X_\cE(X_\cE)^+ \bigl(y - (X_\cE\T)^+ \lambda s \bigr),
\end{equation}
and any lasso solution must be of the form
%
\begin{equation}
\label{eq:lassosol}
\hbeta_{-\cE} = 0\quad \mbox{and}\quad
\hbeta_\cE= (X_\cE)^+ \bigl(y - (X_\cE\T)^+ \lambda s \bigr) + b,
\end{equation}
where $b \in\nul(X_\cE)$. In the case that
$\nul(X_\cE)=\{0\}$---for example, this holds if
$\rank(X)=p$---the lasso solution is unique and is given by
\eqref{eq:lassosol} with $b=0$.
But in general, when $\nul(X_\cE)\not=\{0\}$, it is important to note
that not every $b \in\nul(X_\cE)$ necessarily leads to a lasso
solution in \eqref{eq:lassosol}; the vector $b$ must also preserve
the signs of the nonzero coefficients; that is, it must also satisfy
%
\begin{eqnarray}
\label{eq:lassosign}
\sign\bigl( \bigl[(X_\cE)^+ \bigl(y - (X_\cE\T)^+ \lambda s \bigr) \bigr]_i
+ b_i \bigr) = s_i\nonumber
\\[-8pt]
\\[-8pt]
\eqntext{\mbox{for each $i$ such that }
\bigl[(X_\cE)^+ \bigl(y - (X_\cE\T)^+ \lambda s \bigr) \bigr]_i + b_i
\not= 0.}
\end{eqnarray}
Otherwise, $\gamma$ would not be a proper subgradient of
$\|\hbeta\|_\lone$.

\subsection{Degrees of freedom in terms of the equicorrelation set}
\label{sec:dfequi}

Using relatively simple arguments, we can derive a result on the lasso
degrees of freedom in terms of the equicorrelation set. Our arguments
build on the following key lemma.

\begin{lemma}\label{lem:lassosol0}
For any $y,X$ and $\lambda\geq0$, a lasso solution is given by
%
%
%
\begin{equation}
\label{eq:lassosol0}\hbeta_{-\cE} = 0 \quad\mbox{and}\quad
\hbeta_\cE= (X_\cE)^+ \bigl(y - (X_\cE\T)^+ \lambda s \bigr),
\end{equation}
where $\cE$ and $s$ are the equicorrelation set and signs, as defined in~\eqref{eq:equiset} and~\eqref{eq:equisigns}.
\end{lemma}

In other words, Lemma~\ref{lem:lassosol0} says that the sign condition
\eqref{eq:lassosign} is always satisfied by taking $b=0$,
regardless of the rank of $X$. This result is inspired by the LARS
work of \citet{lars}, though it is not proved in the LARS paper;
see Appendix~B of \citet{ryanphd} for a proof.

Next we show that, almost everywhere in $y$, the equicorrelation set
and signs are locally constant functions of $y$. To emphasize their
functional dependence on $y$, we write them as $\cE(y)$ and $s(y)$.

\begin{lemma}
\label{lem:lcequi}
For almost every $y \in\R^n$, there exists a neighborhood $U$ of $y$
such that $\cE(y')=\cE(y)$ and $s(y)=s(y')$ for all $y' \in U$.
\end{lemma}
\begin{pf}
Define
\[
\cN= \bigcup_{\cE,s} \bigcup_{i \in\cE} \bigl\{z \in\R^n \dvtx
[(X_\cE)^+]_{(i,\cdot)} \bigl(z - (X_\cE\T)^+ \lambda s \bigr) = 0 \bigr\},
\]
where the first union above is taken over all subsets $\cE\subseteq
\{1,\ldots, p\}$ and sign vectors $s \in\{-1,1\}^{|\cE|}$, but we
exclude sets $\cE$ for which a row of $(X_\cE)^+$ is entirely
zero. The set $\cN$ is a finite union of affine subspaces of
dimension $n-1$, and therefore has measure zero.

Let $y \notin\cN$, and abbreviate the equicorrelation set and
signs as $\cE=\cE(y)$ and $s=s(y)$. We may assume no row of
$(X_\cE)^+$ is entirely zero. (Otherwise, this implies that
$X_\cE$ has a zero column, which implies that $\lambda=0$, a
trivial case for this lemma.) Therefore, as $y \notin\cN$, this
means that the lasso solution given in \eqref{eq:lassosol0} satisfies
$\hbeta_i(y) \not= 0$ for every $i \in\cE$.

Now, for a new point $y'$, consider defining
\[
\hbeta_{-\cE}(y') = 0 \quad\mbox{and}\quad
\hbeta_{\cE}(y') = (X_\cE)^+ \bigl(y' - (X_\cE\T)^+ \lambda s \bigr).
\]
We need to verify that this is indeed a solution at $y'$, and
that the corresponding fit has equicorrelation set $\cE$ and signs
$s$. First notice that, after a~straightforward calculation,
\[
X_\cE\T\bigl(y'-X\hbeta(y')\bigr) =
X_\cE\T\bigl(y'- X_\cE(X_\cE)^+ \bigl(y' -
(X_\cE\T)^+\lambda s \bigr) \bigr) = \lambda s.
\]
Also, by the continuity of the function $f \dvtx \R^n \rightarrow\R
^{p-|\cE|}$,
\[
f(x) = X_{-\cE}\T \bigl(x- X_\cE(X_\cE)^+ \bigl(x -
(X_\cE\T)^+\lambda s \bigr) \bigr),
\]
there exists a neighborhood $U_1$ of $y$ such that
\[
\bigl\| X_{-\cE}\T\bigl(y' - X\hbeta(y')\bigr) \bigr\|_\linf
= \bigl\|X_{-\cE}\T \bigl(y'- X_\cE(X_\cE)^+ \bigl(y' -
(X_\cE\T)^+ \lambda s \bigr) \bigr)\bigr \|_\linf< \lambda
\]
for all $y' \in U_1$. Hence $X\hbeta(y')$ has equicorrelation set
$\cE(y')=\cE$ and signs \mbox{$s(y')=s$}.

To check that $\hbeta(y')$ is a lasso solution at $y'$, we consider
the function $g \dvtx \R^n \rightarrow\R^{|\cE|}$,
\[
g(x) = (X_\cE)^+ \bigl(x - (X_\cE\T)^+ \lambda s \bigr).
\]
The continuity of $g$ implies that there exists a neighborhood $U_2$
of $y$ such that
\begin{eqnarray*}
\hbeta_i(y')
&=&
\bigl[(X_\cE\T)^+ \bigl(y' - (X_\cE\T)^+ \lambda s \bigr) \bigr]_i \not= 0
\qquad\mbox{for } i \in\cE, \quad\mbox{and}
\\
\sign(\hbeta_\cE(y'))
&=&
\sign\bigl((X_\cE)^+ \bigl(y' - (X_\cE\T)^+ \lambda s \bigr) \bigr)
\end{eqnarray*}
for each $y' \in U_2$. Defining $U = U_1 \cap U_2$ completes the
proof.
\end{pf}

This immediately implies the following theorem.

\begin{theorem}[(Lasso degrees of freedom, equicorrelation set
representation)]
\label{thm:lassodfequi}
Assume that $y$ follows a normal distribution
\eqref{eq:normal}. For any $X$ and $\lambda\geq0$, the
lasso fit $X\hbeta$ has degrees of freedom
\[
\df(X\hbeta) = \E[\rank(X_\cE)],
\]
where $\cE=\cE(y)$ is the equicorrelation set of the lasso fit at
$y$.
\end{theorem}
\begin{pf}
By Lemmas~\ref{lem:nonexp} and~\ref{lem:lassoproj} we know that
$X\hbeta(y)$ is continuous and almost differentiable, so we can
use Stein's formula \eqref{eq:steindf} for degrees of freedom. By
Lemma~\ref{lem:lcequi}, we know that $\cE=\cE(y)$ and $s=s(y)$ are
locally constant for all $y \notin\cN$. Therefore, taking the
divergence of the fit in \eqref{eq:lassofit}, we get
\[
(\nabla\cdot X\hbeta)(y) =
\tr(X_\cE(X_\cE)^+) =
\rank(X_\cE).
\]
Taking an expectation over $y$ (and recalling that $\cN$ has measure
zero) gives the result.
\end{pf}

Next, we shift our focus to a different subset of variables: the
active set~$\cA$. Unlike the equicorrelation set, the active set is
not unique, as it depends on a~particular choice of lasso
solution. Though it may seem that such nonuniqueness could present
complications, it turns out that all of the active sets share a
special property; namely, the linear subspace $\col(X_\cA)$ is the
same for any choice of active set $\cA$, almost everywhere in
$y$. This invariance allows us to express the degrees of freedom of
lasso fit in terms of the active set (or, more precisely, any active
set).

\subsection{The active set}
\label{sec:lassoact}

Given a particular solution $\hbeta$, we define the \textit{active set}
$\cA$ as
%
\begin{equation}
\cA= \bigl\{ i \in\{1,\ldots, p\} \dvtx \hbeta_i \not= 0
\bigr\}.
\end{equation}
This is also called the support of $\hbeta$ and written $\cA
= \supp(\hbeta)$. From \eqref{eq:lassosol}, we can see that we always
have $\cA\subseteq\cE$, and different active sets $\cA$ can be
formed by choosing $b \in\nul(X_\cE)$ to satisfy the sign condition
\eqref{eq:lassosign} and also
\[
\bigl[(X_\cE)^+ \bigl(y - (X_\cE\T)^+ \lambda s \bigr) \bigr]_i + b_i = 0
\qquad\mbox{for } i \notin\cA.
\]
If $\rank(X)=p$, then $b=0$, so there is a unique active
set, and furthermore $\cA=\cE$ for almost every $y \in\R^n$
(in particular, this last statement holds for \mbox{$y \notin\cN$}, where
$\cN$ is the set of measure zero set defined in the proof of
Lemma~\ref{lem:lcequi}).
For the signs of the coefficients of active variables, we write
%
\begin{equation}
r = \sign(\hbeta_\cA),
\end{equation}
and we note that $r=s_\cA$.

By similar arguments as those used to derive expression
\eqref{eq:lassofit} for the fit in Section~\ref{sec:lassoequi},
the lasso fit can also be written as
%
\begin{equation}
\label{eq:lassofit2}
X\hbeta= (X_\cA)(X_\cA)^+
\bigl(y - (X_\cA\T)^+ \lambda r \bigr)
\end{equation}
for the active set $\cA$ and signs $r$ of any lasso solution
$\hbeta$. If we could take the divergence of the fit in the expression
above, and simply ignore the dependence of $\cA$ and $r$ on $y$ (treat
them as constants), then this would give $(\nabla\cdot X\hbeta)(y) =
\rank(X_\cA)$. In the next section, we show that treating~$\cA$
and~$r$ as constants in \eqref{eq:lassofit2} is indeed
correct, for almost every $y$. This property then implies that the
linear subspace $\col(X_\cA)$ is invariant under any choice of active
set $\cA$, almost everywhere in $y$; moreover, it implies that we can
write the lasso degrees of freedom in terms of any active set.

\subsection{Degrees of freedom in terms of the active set}
\label{sec:dfact}

We first establish a~result on the local stability of $\cA(y)$ and
$r(y)$ [written in this way to emphasize\vspace*{1pt} their dependence on $y$,
through a solution $\hbeta(y)$].

\begin{lemma}
\label{lem:lcact}
There is a set $\cM\subseteq\R^n$, of measure zero, with the
following property:
for $y \notin\cM$, and for any lasso solution $\hbeta(y)$
with active set $\cA(y)$ and signs $r(y)$, there
is a neighborhood $U$ of $y$ such that every point $y' \in U$ yields
a lasso solution $\hbeta(y')$ with the same active set
$\cA(y')=\cA(y)$ and the same active signs $r(y')=r(y)$.
\end{lemma}

The proof is similar to that of Lemma~\ref{lem:lcequi}, except that it
is longer
and somewhat more complicated, so it is delayed until
Appendix~\ref{app:lcact}. Combined with expression \eqref{eq:lassofit2} for
the lasso fit, Lemma~\ref{lem:lcact} now implies an invariance of the subspace spanned by
the active variables.

\begin{lemma}
\label{lem:invact}
For the same set $\cM\subseteq\R^n$ as in Lemma~\ref{lem:lcact}, and
for any $y \notin\cM$, the linear subspace $\col(X_\cA)$ is invariant
under all sets $\cA=\cA(y)$ defined in terms of a lasso solution
$\hbeta(y)$ at $y$.
\end{lemma}
\begin{pf}
Let $y \notin\cM$, and let $\hbeta(y)$ be a solution with active set
$\cA=\cA(y)$ and signs $r=r(y)$. Let $U$ be the neighborhood of $y$ as
constructed in the proof of Lemma~\ref{lem:lcact}; on this
neighborhood, solutions exist with active set $\cA$ and signs
$r$. Hence, recalling \eqref{eq:lassofit2}, we know that for every $y'
\in U$,
\[
X\hbeta(y') = (X_\cA)(X_\cA)^+
\bigl(y' - (X_\cA\T)^+ \lambda r \bigr).
\]
Now suppose that $\cA^*$ and $r^*$ are the active set and signs of
another lasso solution at $y$. Then, by the same arguments, there is
a neighborhood $U^*$ of~$y$ such that
\[
X\hbeta(y') = (X_{\cA^*})(X_{\cA^*})^+
\bigl(y' - (X_{\cA^*}\T)^+ \lambda r^* \bigr)
\]
for all $y' \in U^*$.
By the uniqueness of the fit, we have that for each $y' \in U \cap
U^*$,
\[
(X_\cA)(X_\cA)^+
\bigl(y' - (X_\cA\T)^+ \lambda r \bigr) =
(X_{\cA^*})(X_{\cA^*})^+
\bigl(y' - (X_{\cA^*}\T)^+ \lambda r^* \bigr).
\]
Since $U \cap U^*$ is open, for any $z \in\col(X_\cA)$, there
is an $\varepsilon>0$ such that $y+\varepsilon z \in U \cap U^*$. Plugging
$y'=y+\varepsilon z$ into the above equation implies that $z \in
\col(X_{\cA^*})$, so $\col(X_\cA) \subseteq\col(X_{\cA^*})$. A
similar argument gives $\col(X_{\cA^*}) \subseteq\col(X_\cA)$,
completing the proof.
\end{pf}

Again, this immediately leads to the following theorem.

\begin{theorem}[(Lasso degrees of freedom, active set
representation)]
\label{thm:lassodfact}
Assume that $y$ follows a normal distribution
\eqref{eq:normal}. For any $X$ and $\lambda\geq0$, the lasso
fit $X\hbeta$ has degrees of freedom
\[
\df(X\hbeta) = \E[\rank(X_\cA)],
\]
where $\cA=\cA(y)$ is the active set corresponding to any lasso
solution $\hbeta(y)$ at~$y$.
\end{theorem}

Note: By Lemma~\ref{lem:invact}, $\rank(X_\cA)$ is an
invariant quantity, not depending on the choice of
active set (coming from a lasso solution), for almost every
$y$. This makes the above result well defined.
\begin{pf*}{Proof of Theorem~\ref{thm:lassodfact}}
We can apply Stein's formula \eqref{eq:steindf} for degrees of
freedom, because $X\hbeta(y)$ is continuous and almost differentiable
by Lemmas~\ref{lem:nonexp} and~\ref{lem:lassoproj}. Let $\cA=\cA(y)$
and $r=r(y)$ be the active set and active signs of a lasso solution
at $y \notin\cM$, with $\cM$ as in Lemma~\ref{lem:invact}. By this
same lemma, there exists a lasso solution with active set $\cA$ and
signs $r$ at every point $y'$ in some neighborhood $U$ of $y$, and
therefore, taking the divergence of the fit~\eqref{eq:lassofit2}, we
get
\[
(\nabla\cdot X\hbeta)(y) = \tr(X_\cA(X_\cA)^+) = \rank(X_\cA).
\]
Taking an expectation over $y$ completes the proof.
\end{pf*}

\begin{remark*}[(Equicorrelation set representation)]
The proof of Lemma~\ref{lem:lcact} showed that, for almost every $y$, the equicorrelation
set $\cE$ is actually the active set $\cA$ of the particular lasso
solution defined in \eqref{eq:lassosol0}. Hence Theorem~\ref{thm:lassodfequi}
can be viewed as a corollary of Theorem~\ref{thm:lassodfact}.
\end{remark*}

\begin{remark*}[(Full column rank $X$)]
When $\rank(X)=p$, the lasso solution is unique, and there is only one
active set $\cA$. And as the columns of $X$ are linearly independent,
we have $\rank(X)=|\cA|$, so the result of Theorem~\ref{thm:lassodfact} reduces to
\[
\df(X\hbeta) = \E|\cA|,
\]
as shown in \citet{lassodf}.
\end{remark*}

\begin{remark*}[(The smallest active set)]
 An interesting result on the lasso
degrees of freedom was recently and independently obtained by \citet{dfl1}.
Their result states that, for a general $X$,
\[
\df(X\hbeta) = \E|\cA^*|,
\]
where $|\cA^*|$ is the smallest cardinality among all active
sets of lasso solutions. This actually follows from Theorem~\ref
{thm:lassodfact},
by noting that for any $y$ there exists a lasso solution\vadjust{\goodbreak} whose
active set $\cA^*$ corresponds to linear independent predictors
$X_{\cA^*}$,
so $\rank(X_{\cA^*})=|\cA^*|$ [e.g., see Theorem 3 in Appendix~B of
\citet{boostpath}], and furthermore, for almost every $y$ no active set
can have a cardinality smaller than $|\cA^*|$, as this would
contradict Lemma~\ref{lem:invact}.
\end{remark*}

\begin{remark*}[(The elastic net)]
 Consider the elastic net problem
[\citet{enet}],
%
\begin{equation}
\label{eq:enet}
\hbeta= \mathop\argmin_{\beta\in\R^p} \half\|y-X\beta\|_\ltwo^2 +
\lambda_1 \|\beta\|_\lone+ \frac{\lambda_2}{2} \|\beta\|_\ltwo^2,
\end{equation}
where we now have two tuning parameters $\lambda_1,\lambda_2 \geq
0$. Note that our notation above emphasizes the fact that there is
always a unique solution to the elastic net criterion, regardless of
the rank of $X$. This property (among others, such as stability and
predictive ability) is considered an advantage of the elastic net over
the lasso. We can rewrite the elastic net problem
\eqref{eq:enet} as a (full column rank) lasso problem,
\[
\hbeta= \mathop\argmin_{\beta\in\R^p} \half \left\|
\left(
\matrix{ y
\cr 0}
\right) -
\left[
\matrix{X
\cr
 \sqrt{\lambda_2} I}
\right]
\beta \right\|_\ltwo^2 + \lambda_1 \|\beta\|_\lone,
\]
and hence it can be shown (although we omit the details)
that the degrees of freedom of the elastic net fit is
\[
\df(X\hbeta) = \E\bigl[\tr\bigl(X_\cA(X_\cA\T X_\cA+ \lambda_2
I)^{-1} X_\cA\T\bigr) \bigr],
\]
where $\cA=\cA(y)$ is the active set of the elastic net solution at
$y$.
\end{remark*}

\begin{remark*}[(The lasso with intercept)]
 It is often more appropriate
to include an (unpenalized) intercept coefficient in the lasso model,
yielding the problem
%
\begin{equation}
\label{eq:lassoint}
(\hbeta_0,\hbeta) \in\mathop\argmin_{(\beta_0,\beta) \in\R^{p+1}}
\half
\|y-\beta_0\mathbh{1} -X\beta\|_\ltwo^2 + \lambda\|\beta\|_\lone,
\end{equation}
where $\mathbh{1}=(1,1,\ldots,1) \in\R^n$ is the vector of all
$1$s. Defining $M=I-\mathbh{1}\mathbh{1}\T/n \in\R^{n \times n}$, we
note that the fit of problem \eqref{eq:lassoint} can be written as
$\hbeta_0\mathbh{1} + X\hbeta= (I-M)y + MX\hbeta$, and that
$\hbeta$ solves the usual lasso problem
\[
\hbeta\in\mathop\argmin_{\beta\in\R^p} \half
\|My - MX\beta\|_\ltwo^2 + \lambda\|\beta\|_\lone.
\]
Now it follows (again we omit the details) that the fit of the
lasso problem with intercept \eqref{eq:lassoint} has degrees of
freedom
\[
\df(\hbeta_0 \mathbh{1} + X\hbeta) = 1 + \E[\rank(MX_\cA)],
\]
where $\cA=\cA(y)$ is the active set of a solution $\hbeta(y)$ at $y$
(these are the nonintercept coefficients). In other words, the
degrees of freedom is one plus the expected\vadjust{\goodbreak} dimension of the subspace
spanned by the active variables, once we have centered these
variables. A similar result holds for an arbitrary set of unpenalized
coefficients, by replacing $M$ above with the projection onto the
orthogonal complement of the column space of the unpenalized
variables, and $1$ above with the dimension of the column space of the
unpenalized variables.
\end{remark*}

As mentioned in the \hyperref[sec:intro]{Introduction}, a nice feature of the full column
rank result \eqref{eq:lassodffull} is its interpretability and its
explicit nature. The general result is also explicit in the sense that
an unbiased estimate of degrees of freedom can be achieved by computing
the rank of a given matrix. In terms of interpretability, when
$\rank(X)=p$, the degrees of
freedom of the lasso fit is $\E|\cA|$---this says that, on average,
the lasso ``spends'' the same number of parameters as does
linear regression on $|\cA|$ linearly independent predictor
variables. Fortunately, a similar interpretation is possible in the
general case: we showed in Theorem~\ref{thm:lassodfact} that for a
general predictor matrix $X$, the degrees of freedom of the lasso fit
is $\E[\rank(X_\cA)]$, the expected dimension of the linear subspace
spanned by the active variables. Meanwhile, for the linear regression
problem
%
\begin{equation}
\label{eq:lsa}
\hbeta_\cA= \mathop\argmin_{\beta_\cA\in\R^{|\cA|}}
\|y-X_\cA\beta_\cA\|_\ltwo^2,
\end{equation}
where we consider $\cA$ fixed, the degrees of freedom of the fit is
$\tr(X_\cA(X_\cA)^+) = \rank(X_\cA)$. In
other words, the lasso adaptively selects a subset $\cA$ of the
variables to use for a linear model of $y$, but on average it only
``spends'' the same number of parameters as would linear regression on
the variables in $\cA$, if $\cA$ was pre-specified.

How is this possible? Broadly speaking, the answer lies in the
shrinkage due to the $\ell_1$ penalty. Although the active set is
chosen adaptively, the lasso does not estimate the active coefficients
as aggressively as does the corresponding linear regression problem
\eqref{eq:lsa}; instead, they are shrunken toward zero, and this
adjusts for the adaptive selection.
Differing views have been presented in the literature with respect to
this feature of lasso shrinkage. On the one hand, for example,
\citet{scad} point out that lasso estimates suffer from bias due to the
shrinkage of large coefficients, and motivate the nonconvex \textit{SCAD}
penalty as an attempt to overcome this bias. On the other hand, for
example, Loubes and Massart (\citeyear{loubes}) discuss the merits of such shrunken
estimates in model selection criteria, such as~\eqref{eq:tunsel}. In
the current context, the shrinkage due to the $\ell_1$ penalty is helpful
in that it provides control over degrees of freedom. A more precise
study of
this idea is the topic of future work.

\section{The generalized lasso}
\label{sec:genlasso}

In this section we extend our degrees of freedom results to the
generalized lasso problem, with an arbitrary predictor matrix $X$ and
penalty matrix $D$. As before, the KKT conditions play a~central role,\vadjust{\goodbreak} and we present these first. Also, many results that follow
have equivalent derivations from the perspective of the generalized lasso
dual problem; see Appendix~\ref{app:dual}.
We remind the reader that $D_R$ is used to extract to extract rows of
$D$ corresponding to an index set $R$.

\subsection{The KKT conditions and the underlying polyhedron}
\label{sec:genlassokkt}

The KKT conditions for the generalized lasso problem
\eqref{eq:genlasso} are
%
\begin{eqnarray}
\label{eq:genlassokkt}
&\displaystyle X\T( y-X\hbeta) = D\T\lambda\gamma,&
\\
\label{eq:genlassosg}
&\displaystyle\gamma_i \in
\cases{
\{\sign((D\hbeta)_i)\} & \quad if $(D\hbeta)_i \not= 0$,
 \cr
[-1,1] & \quad  if $(D\hbeta)_i = 0$.}&
\end{eqnarray}
Now $\gamma\in\R^m$ is a subgradient of the function
$f(x)=\|x\|_\lone$ evaluated at $x=D\hbeta$. Similar to what we
showed for the lasso, it follows from the KKT conditions that the
generalized lasso fit is the residual from projecting $y$ onto a
polyhedron.

\begin{lemma}
\label{lem:genlassoproj}
For any $X$ and $\lambda\geq0$, the generalized lasso fit
can be written as $X\hbeta(y) = (I-P_C)(y)$, where $C
\subseteq\R^n$ is the polyhedron
\[
C = \{ u \in\R^n \dvtx X\T u = D\T w \mbox{ for }
w \in\R^m, \|w\|_\linf\leq\lambda\}.
\]
\end{lemma}
\begin{pf}
The proof is quite similar to that of Lemma~\ref{lem:lassoproj}. As in \eqref{eq:innerprod}, we want to show that
%
\begin{equation}
\label{eq:innerprod2}
\langle X\hbeta, y-X\hbeta\rangle
- \langle X\T u, \hbeta\rangle\geq0
\end{equation}
for all $u \in C$, where $C$ is as in the lemma. For the first
term above, we can take an inner product with $\hbeta$ on both sides
of \eqref{eq:genlassokkt} to get $\langle X\hbeta, y-X\hbeta\rangle=
\lambda\|D\hbeta\|_\lone$, and furthermore,
\[
\lambda\|D\hbeta\|_\lone= \max_{\|w\|_\linf\leq\lambda} \langle w,
D\hbeta\rangle= \max_{\|w\|_\linf\leq\lambda} \langle D\T w,
\hbeta\rangle.
\]
Therefore \eqref{eq:innerprod2} holds if $X\T u = D\T w$
for some $\|w\|_\linf\leq\lambda$, in other words, if $u
\in C$. To show that $C$ is a polyhedron, note that we can
write it as $C = (X\T)^{-1} (D\T(B))$
where $(X\T)^{-1}$ is taken to mean the inverse image
under the linear map $X\T$,
and $B=\{ w \in\R^m \dvtx \|w\|_\linf\leq\lambda\}$, a hypercube in
$\R^m$. Clearly~$B$ is a polyhedron, and
the image or inverse image of a~polyhedron under a~linear
map is still a polyhedron.
\end{pf}

As with the lasso, this lemma implies that the generalized
lasso fit $X\hbeta(y)$ is nonexpansive, and therefore continuous and
almost differentiable as a~function of $y$, by Lemma~\ref{lem:nonexp}.
This is important because it allows us to use Stein's formula
when computing degrees of freedom.

In the next section we define the boundary set $\cB$, and derive
expressions for the generalized lasso fit and solutions in terms of
$\cB$. The following section defines the active set $\cA$ in the
generalized lasso context, and again gives expressions for the fit\vadjust{\goodbreak}
and solutions in terms of $\cA$. Though neither $\cB$ nor~$\cA$ are
necessarily unique for the generalized lasso problem, any choice of~$\cB$ or~$\cA$ generates a special invariant subspace (similar to the
case for the active sets in the lasso problem). We are subsequently
able to express the degrees of freedom of the generalized lasso fit in
terms of any boundary set~$\cB$, or any active set $\cA$.

\subsection{The boundary set}
\label{sec:genlassobound}

Like the lasso, the generalized lasso fit $X\hbeta$ is always unique
(following from Lemma~\ref{lem:genlassoproj}, and the fact that
projection onto a~closed convex set is unique).
However, unlike the lasso, the optimal subgradient $\gamma$
in the generalized lasso problem is not
necessarily unique. In particular, if $\rank(D)<m$, then the
optimal subgradient $\gamma$ is not uniquely determined by
conditions \eqref{eq:genlassokkt} and \eqref{eq:genlassosg}.
Given a subgradient $\gamma$ satisfying~\eqref{eq:genlassokkt} and \eqref{eq:genlassosg} for some $\hbeta$,
we define the \textit{boundary set} $\cB$ as
\[
\cB= \{ i \in\{1,\ldots, m\} \dvtx |\gamma_i|=1 \}.
\]
This generalizes the notion of the equicorrelation set $\cE$ in the
lasso problem [though, as just noted, the set $\cB$ is not necessarily
unique unless \mbox{$\rank(D)=m$}]. We also define
\[
s = \gamma_\cB.
\]

Now we focus on writing the generalized lasso fit and solutions
in terms of~$\cB$ and~$s$. Abbreviating $P=P_{\nul(D_{-\cB})}$, note
that we can expand $P D\T\lambda\gamma= P D_\cB\T\lambda s + P
D_{-\cB}\T\lambda\gamma_{-\cB} = P D_\cB\T\lambda s$. Therefore,
multiplying both sides of \eqref{eq:genlassokkt} by $P$ yields
%
\begin{equation}
\label{eq:genlassokktb}
P X\T(y-X\hbeta) = P D_\cB\T\lambda s.
\end{equation}
Since $P D_\cB\T\lambda s \in\col(PX\T)$, we can write\vspace*{1pt}
$P D_\cB\T\lambda s = (P X\T)(P X\T)^+ P D_\cB\T\lambda s
= (P X\T)(P X\T)^+ D_\cB\T\lambda s$. Also, we have
$D_{-\cB}\hbeta=0$ by definition of $\cB$, so $P \hbeta=
\hbeta$. These two facts allow us to rewrite
\eqref{eq:genlassokktb} as
\[
P X\T X P \hbeta= P X\T
\bigl(y - (P X\T)^+ D_\cB\T\lambda s \bigr),
\]
and hence the fit $X\hbeta= X P \hbeta$ is
%
\begin{equation}
\label{eq:genlassofit}
X\hbeta= \bigl(X P_{\nul(D_{-\cB})}\bigr)\bigl(X P_{\nul(D_{-\cB})}\bigr)^+
\bigl(y - \bigl(P_{\nul(D_{-\cB})} X\T\bigr)^+ D_\cB\T\lambda s \bigr),
\end{equation}
where we have un-abbreviated $P=P_{\nul(D_{-\cB})}$. Further, any
generalized lasso solution is of the form
%
\begin{equation}
\label{eq:genlassosol}
\hbeta= \bigl(X P_{\nul(D_{-\cB})}\bigr)^+
\bigl(y - \bigl(P_{\nul(D_{-\cB})} X\T\bigr)^+ D_\cB\T\lambda s \bigr) + b,
\end{equation}
where $b \in\nul(X P_{\nul(D_{-\cB})})$. Multiplying the above
equation by $D_{-\cB}$, and recalling that $D_{-\cB}\hbeta=0$,
reveals that $b \in\nul(D_{-\cB})$; hence
$b \in\nul(XP_{\nul(D_{-\cB})}) \cap\nul(D_{-\cB}) =
\nul(X) \cap\nul(D_{-\cB})$. In the case that $\nul(X)
\cap\nul(D_{-\cB}) = \{0\}$, the generalized lasso solution is
unique and is given by \eqref{eq:genlassosol} with $b=0$. This
occurs when $\rank(X)=p$, for example. Otherwise, any $b \in\nul(X)
\cap\nul(D_{-\cB})$ gives\vadjust{\goodbreak} a generalized lasso solution in
\eqref{eq:genlassosol} as long as it also satisfies the sign
condition
%
\begin{eqnarray}
\label{eq:genlassosign}
&&\hspace*{37pt}\sign\bigl(D_i \bigl(X P_{\nul(D_{-\cB})}\bigr)^+
\bigl(y - \bigl(P_{\nul(D_{-\cB})} X\T\bigr)^+ D_\cB\T\lambda s \bigr) +
D_i b \bigr) = s_i\nonumber
\\
&&\hspace*{43pt}\mbox{for each } i \in\cB \mbox{ such that }
D_i \bigl(X P_{\nul(D_{-\cB})}\bigr)^+
\\
&&\hspace*{139.7pt}\qquad{}\times\bigl(y - \bigl(P_{\nul(D_{-\cB})} X\T\bigr)^+ D_\cB\T\lambda s \bigr) +
D_i b \not= 0,\nonumber
\end{eqnarray}
necessary to ensure that $\gamma$ is a proper subgradient of
$\|D\hbeta\|_\lone$.

\subsection{The active set}
\label{sec:genlassoact}

We define the \textit{active set} of a particular solution~$\hbeta$~as
\[
\cA= \{ i \in\{1,\ldots, m\} \dvtx (D\hbeta)_i \not= 0 \},
\]
which can be alternatively expressed as $\cA=\supp(D\hbeta)$.
If $\hbeta$ corresponds to a subgradient with
boundary set $\cB$ and signs $s$, then $\cA\subseteq\cB$;
in particular, given $\cB$ and $s$, different active sets $\cA$ can be
generated by taking $b \in\nul(X)\cap\nul(D_{-\cB})$ such that
\eqref{eq:genlassosign} is satisfied, and also
\[
D_i \bigl(X P_{\nul(D_{-\cB})}\bigr)^+
\bigl(y - \bigl(P_{\nul(D_{-\cB})} X\T\bigr)^+ D_\cB\T\lambda s \bigr) +
D_i b = 0 \qquad\mbox{for } i \in\cB\setminus\cA.
\]
If $\rank(X)=p$, then $b=0$, and there is only one active set
$\cA$; however, in this case, $\cA$ can still be a strict
subset of $\cB$. This is quite different from the lasso problem,
wherein $\cA=\cE$ for almost every $y$ whenever $\rank(X)=p$.
[Note that in the generalized lasso problem, $\rank(X)=p$ implies that
$\cA$ is unique but implies nothing about the uniqueness of
$\cB$---this is determined by the rank of $D$. The boundary set $\cB$
is not necessarily unique if $\rank(D)<m$, and in this case we may
have $D_i (X P_{\nul(D_{-\cB})})^+=0$ for some $i \in\cB$, which
certainly implies that $i \notin\cA$ for any $y \in\R^n$. Hence some
boundary sets may not correspond to active sets at any~$y$.]
We denote the signs of the
active entries in $D\hbeta$ by
\[
r = \sign(D_\cA\hbeta),
\]
and we note that $r=s_\cA$.

Following the same arguments as those leading up to the expression for
the
fit \eqref{eq:genlassofit} in Section~\ref{sec:genlassobound},
we can alternatively express the generalized lasso fit as
%
\begin{equation}
\label{eq:genlassofit2}\qquad
X\hbeta= \bigl(X P_{\nul(D_{-\cA})}\bigr)\bigl(X P_{\nul(D_{-\cA})}\bigr)^+
\bigl(y - \bigl(P_{\nul(D_{-\cA})} X\T\bigr)^+ D_\cA\T\lambda r \bigr),
\end{equation}
where $\cA$ and $r$ are the active set and signs of any solution.
Computing the divergence of the fit in \eqref{eq:genlassofit2},
and pretending that $\cA$ and $r$ are constants (not depending on
$y$), gives $(\nabla\cdot X\hbeta)(y) =
\dim(\col(XP_{\nul(D_{-\cA})})) = \dim(X(\nul(D_{-\cA})))$. The same
logic applied to \eqref{eq:genlassofit} gives
$(\nabla\cdot X\hbeta)(y) = \dim(X(\nul(D_{-\cB})))$. The next
section shows that, for almost every $y$, the quantities $\cA,r$ or
$\cB,s$ can indeed be treated as locally constant in expressions
\eqref{eq:genlassofit2} or \eqref{eq:genlassofit}, respectively.
We then prove that linear subspaces
$X(\nul(D_{-\cB})), X(\nul(D_{-\cA}))$ are invariant under all
choices of boundary sets $\cB$, respectively active sets $\cA$, and
that the two subspaces are in fact equal, for almost every $y$.
Furthermore, we express the generalized lasso degrees of freedom
in terms of any boundary set or any active set.

\subsection{Degrees of freedom}
\label{sec:genlassodf}

We call $(\gamma(y),\hbeta(y))$ an \textit{optimal pair} provided
that $\gamma(y)$ and $\hbeta(y)$ jointly satisfy the KKT conditions,
\eqref{eq:genlassokkt} and \eqref{eq:genlassosg}, at $y$. For such a
pair, we consider its boundary set $\cB(y)$, boundary signs $s(y)$,
active set $\cA(y)$, active signs $r(y)$, and show that these sets
and sign vectors possess a kind of local stability.

\begin{lemma}
\label{lem:lcbound}
There exists a set $\cN\subseteq\R^n$, of measure zero, with the
following property: for $y \notin\cN$, and for any
optimal pair $(\gamma(y),\hbeta(y))$ with boundary set
$\cB(y)$, boundary signs $s(y)$, active set $\cA(y)$, and active signs
$r(y)$, there is a neighborhood $U$ of $y$ such that each point $y'
\in U$ yields an optimal pair $(\gamma(y'),\hbeta(y'))$ with the same
boundary set $\cB(y')=\cB(y)$, boundary signs $s(y')=s(y)$,
active set $\cA(y')=\cA(y)$ and active signs $r(y')=r(y)$.
\end{lemma}

The proof is delayed to Appendix~\ref{app:lcbound}, mainly because of
its length. Now Lemma~\ref{lem:lcbound}, used together with
expressions \eqref{eq:genlassofit} and \eqref{eq:genlassofit2} for
the generalized lasso fit, implies an invariance in representing a
(particularly important) linear subspace.

\begin{lemma}
\label{lem:invbound}
For the same set $\cN\subseteq\R^n$ as in Lemma~\ref{lem:lcbound},
and for any $y \notin\cN$, the linear subspace
$L=X(\nul(D_{-\cB}))$ is invariant under all boundary sets
$\cB=\cB(y)$ defined in terms of an optimal subgradient at $\gamma(y)$
at $y$. The linear subspace $L'=X(\nul(D_{-\cA}))$ is also
invariant under all choices of active sets $\cA=\cA(y)$ defined in
terms of a generalized lasso solution $\hbeta(y)$ at~$y$. Finally,
the two subspaces are equal, $L=L'$.
\end{lemma}
\begin{pf}
Let $y \notin\cN$, and let $\gamma(y)$ be an optimal subgradient with
boundary set $\cB=\cB(y)$ and signs $s=s(y)$. Let $U$ be the
neighborhood of $y$ over which optimal subgradients exist with
boundary set $\cB$ and signs $s$, as given by Lemma~\ref{lem:lcbound}. Recalling the expression for the fit
\eqref{eq:genlassofit}, we have that for every $y' \in U$
\[
X\hbeta(y') =
\bigl(X P_{\nul(D_{-\cB})}\bigr)\bigl(X P_{\nul(D_{-\cB})}\bigr)^+
\bigl(y' - \bigl(P_{\nul(D_{-\cB})} X\T\bigr)^+ D_\cB\T\lambda s \bigr).
\]
If $\hbeta(y)$ is a solution with active set $\cA=\cA(y)$ and
signs $r=r(y)$, then again by Lemma~\ref{lem:lcbound} there is a
neighborhood $V$ of $y$ such that each point $y' \in V$ yields a
solution with active set $\cA$ and signs $r$. [Note that $V$
and $U$ are not necessarily equal unless $\gamma(y)$ and $\hbeta(y)$
jointly satisfy the KKT conditions at~$y$.] Therefore, recalling
\eqref{eq:genlassofit}, we have
\[
X\hbeta(y') =
\bigl(X P_{\nul(D_{-\cA})}\bigr)\bigl(X P_{\nul(D_{-\cA})}\bigr)^+
\bigl(y' - \bigl(P_{\nul(D_{-\cA})} X\T\bigr)^+ D_\cA\T\lambda r \bigr)
\]
for each $y' \in V$. The uniqueness of the generalized lasso fit now
implies that
\begin{eqnarray*}
&&\bigl(X P_{\nul(D_{-\cB})}\bigr)\bigl(X P_{\nul(D_{-\cB})}\bigr)^+
\bigl(y' - \bigl(P_{\nul(D_{-\cB})} X\T\bigr)^+ D_\cB\T\lambda s \bigr)
\\
&&\qquad=
\bigl(X P_{\nul(D_{-\cA})}\bigr)\bigl(X P_{\nul(D_{-\cA})}\bigr)^+
\bigl(y' - \bigl(P_{\nul(D_{-\cA})} X\T\bigr)^+ D_\cA\T\lambda r \bigr)
\end{eqnarray*}
for all $y' \in U \cap V$.
As $U \cap V$ is open, for any $z \in\col(XP_{\nul(D_{-\cB})})$,
there exists an $\varepsilon>0$ such that $y+\varepsilon z \in U \cap
V$. Plugging $y'=y+\varepsilon z$ into the equation above reveals
that $z \in\col(XP_{\nul(D_{-\cA})})$, hence
$\col(XP_{\nul(D_{-\cB})}) \subseteq\col(XP_{\nul(D_{-\cA})})$. The
reverse inclusion follows similarly, and therefore\break
$\col(XP_{\nul(D_{-\cB})}) = \col(XP_{\nul(D_{-\cA})})$. Finally, the
same strategy can be used to show that these linear subspaces are
unchanged for any choice of boundary set $\cB=\cB(y)$, coming from an
optimal subgradient at $y$ and for any choice of active
set $\cA=\cA(y)$ coming from a solution at $y$. Noticing
that $\col(M P_{\nul(N)}) = M(\nul(N))$ for matrices $M,N$ gives the
result as stated in the lemma.
\end{pf}

This local stability result implies the following theorem.

\begin{theorem}[(Generalized lasso degrees of freedom)]
\label{thm:genlassodf}
Assume that $y$ follows a normal distribution
\eqref{eq:normal}. For any $X,D$ and $\lambda\geq0$, the degrees of
freedom of the generalized lasso fit can be expressed as
\[
\df(X\hbeta) = \E[\dim(X(\nul(D_{-\cB})))],
\]
where $\cB=\cB(y)$ is the boundary set corresponding to any optimal
subgradient~$\gamma(y)$ of the generalized lasso problem at $y$.
We can alternatively express degrees of freedom as
\[
\df(X\hbeta) = \E[\dim(X(\nul(D_{-\cA})))],
\]
with $\cA=\cA(y)$ being the active set corresponding to any
generalized lasso solution $\hbeta(y)$ at $y$.
\end{theorem}

Note: Lemma~\ref{lem:invbound} implies that for almost every $y \in
\R^n$, for any $\cB$ defined in terms of an optimal subgradient, and
for any $\cA$ defined in terms of a~generalized lasso solution,
$\dim(X(\nul(D_{-\cB})))=\dim(X(\nul(D_{-\cA})))$.
This makes the above expressions for degrees of freedom well defined.
\begin{pf*}{Proof of Theorem~\ref{thm:genlassodf}}
First, the continuity and almost differentiability of $X\hbeta(y)$
follow from Lemmas~\ref{lem:nonexp} and~\ref{lem:genlassoproj}, so we
can use Stein's formula \eqref{eq:steindf} for degrees of freedom.
Let $y \notin\cN$, where $\cN$ is the set of measure zero as in Lemma~\ref{lem:lcact}. If $\cB=\cB(y)$ and $s=s(y)$ are the boundary set and
signs of an optimal subgradient at $y$, then by Lemma~\ref{lem:invbound} there is a
neighborhood $U$ of~$y$ such that each\vadjust{\goodbreak}
point $y' \in U$ yields an optimal subgradient with boundary set $\cB$
and signs $s$. Therefore, taking the divergence of the fit in
\eqref{eq:genlassofit},
\[
(\nabla\cdot X\hbeta)(y) = \tr\bigl(P_{X(\nul(D_{-\cB}))}\bigr) =
\dim(X(\nul(D_{-\cB}))),
\]
and taking an expectation over $y$ gives the first expression in the
theorem.

Similarly, if $\cA=\cA(y)$ and $r=r(y)$ are the active set and signs
of a~generalized lasso solution at $y$, then by Lemma~\ref{lem:invbound} there exists a solution with active set $\cA$ and
signs $r$ at each point $y'$ in some neighborhood $V$ of $y$.
The divergence of the fit in \eqref{eq:genlassofit2} is hence
\[
(\nabla\cdot X\hbeta)(y) = \tr\bigl(P_{X(\nul(D_{-\cA}))}\bigr) =
\dim(X(\nul(D_{-\cA}))),
\]
and taking an expectation over $y$ gives the second expression.
\end{pf*}

\begin{remark*}[(Full column rank $X$)]
If $\rank(X)=p$, then $\dim(X(L))=\dim(L)$ for any linear subspace
$L$, so the results of Theorem~\ref{thm:genlassodf} reduce to
\[
\df(X\hbeta) = \E[\nuli(D_{-\cB})] = \E[\nuli(D_{-\cA})].
\]
The first equality above was shown in \citet{genlasso}.
Analyzing the null space of $D_{-\cB}$ (equivalently, $D_{-\cA}$) for
specific choices of $D$ then gives interpretable results on the
degrees of freedom of the fused lasso and trend filtering fits
as mentioned in the introduction. It is important to note that, as
$\rank(X)=p$, the active set $\cA$ is unique, but not necessarily
equal to the boundary set $\cB$ [since $\cB$ can be nonunique if
$\rank(D)<m$].
\end{remark*}

\begin{remark*}[(The lasso)]
If $D=I$, then $X(\nul(D_{-S})) = \col(X_S)$ for any subset $S
\subseteq\{1,\ldots, p\}$. Therefore the results of Theorem~\ref{thm:genlassodf} become
\[
\df(X\hbeta) = \E[\rank(X_\cB)] = \E[\rank(X_\cA)],
\]
which match the results of Theorems~\ref{thm:lassodfequi} and~\ref{thm:lassodfact} (recall that for the lasso the boundary set $\cB$
is exactly the same as equicorrelation set $\cE$).
\end{remark*}

\begin{remark*}[(The smallest active set)]
Recent and independent work of \citet{analreg} shows that, for
arbitrary $X,D$ and for any $y$, there exists a generalized lasso solution
whose active set $\cA^*$ satisfies
\[
\nul(X) \cap\nul(D_{-\cA^*}) = \{0\}.
\]
(Calling $\cA^*$ the ``smallest'' active set is somewhat of an abuse of
terminology, but it is the smallest in terms of the above intersection.)
The authors then prove that, for any $X,D$, the generalized lasso fit has
degrees of freedom
\[
\df(X\hbeta) = \E[\nuli(D_{-\cA^*})],
\]
with $\cA^*$ the special active set as above. This matches the active set
result of Theorem~\ref{thm:genlassodf} applied to $\cA^*$, since
$\dim(X(\nul(D_{-\cA^*}))) = \nuli(D_{-\cA^*})$ for this special
active set.
\end{remark*}

We conclude this section by comparing the active set result of Theorem~\ref{thm:genlassodf}
to degrees of freedom in a particularly relevant\vadjust{\goodbreak}
equality constrained linear regression problem (this comparison is
similar to that made in lasso case, given at the end of Section~\ref{sec:lasso}). The result states that the generalized lasso fit
has degrees of freedom $\E[\dim(X(\nul(D_{-\cA})))]$, where
$\cA=\cA(y)$ is the active set of a generalized lasso solution at
$y$. In other words, the complement of $\cA$ gives the rows of $D$
that are orthogonal to some generalized lasso solution. Now, consider
the equality constrained linear regression problem
%
\begin{equation}
\label{eq:lsb}
\hbeta\in\mathop\argmin_{\beta\in\R^p}
\|y-X\beta\|_\ltwo^2
\qquad\mbox{subject to } D_{-\cA}\beta=0,
\end{equation}
in which the set $\cA$ is fixed.
It is straightforward to verify that the fit of this problem is
the projection map onto
$\col(XP_{\nul(D_{-\cA})})=X(\nul(D_{-\cA}))$,
and hence has degrees of freedom $\dim(X(\nul(D_{-\cA})))$. This means
that the generalized lasso fits a linear model of $y$, and
simultaneously makes the coefficients orthogonal to an adaptive subset
$\cA$ of the rows of $D$, but on average it only uses the same
number of parameters as does the corresponding equality constrained
linear regression problem \eqref{eq:lsb}, in which $\cA$ is
pre-specified.

This seemingly paradoxical statement can be explained by the shrinkage
due to the $\ell_1$ penalty. Even though the active set $\cA$ is
chosen adaptively based on $y$, the generalized lasso does not
estimate the coefficients as aggressively as does the equality
constrained linear regression problem \eqref{eq:lsb}, but rather, it
shrinks them toward zero. Roughly speaking, his shrinkage can be
viewed as a ``deficit'' in degrees of freedom, which makes up for the
``surplus'' attributed to the adaptive selection. We study this idea
more precisely in a future paper.

\section{Discussion}
\label{sec:disc}

We showed that the degrees of freedom of the lasso fit, for an
arbitrary predictor matrix $X$, is equal to $\E[\rank(X_\cA)]$.
Here $\cA=\cA(y)$ is the active set of any lasso
solution at $y$, that is, $\cA(y)=\supp(\hbeta(y))$. This result is
well defined, since we proved that any active set $\cA$ generates
the same linear subspace $\col(X_\cA)$, almost everywhere in $y$. In
fact, we showed that for almost every $y$, and for any active set
$\cA$ of a solution at $y$, the lasso fit can be written as
\[
X\hbeta(y') = P_{\col(X_\cA)}(y') + c
\]
for all $y'$ in a neighborhood of $y$, where $c$ is a constant
(it does not depend on $y'$). This draws an interesting connection to
linear regression, as it shows that locally the lasso fit is just
a translation of the linear regression fit of on~$X_\cA$. The same
results (on degrees of freedom and local representations of the
fit) hold when the active set $\cA$ is replaced by the
equicorrelation set~$\cE$.

Our results also extend to the generalized lasso
problem, with an arbitrary predictor matrix $X$ and arbitrary penalty
matrix $D$. We showed that degrees of freedom of the generalized lasso
fit is $\E[\dim(X(\nul(D_{-\cA})))]$, with $\cA=\cA(y)$ being the
active set of any generalized lasso solution at $y$, that is,
$\cA(y)=\supp(D\hbeta(y))$. As before, this result is
well defined because any choice of active set $\cA$ generates the
same linear subspace $X(\nul(D_{-\cA}))$, almost everywhere in $y$.
Furthermore, for almost every $y$, and for any active set of a
solution at $y$, the generalized lasso fit satisfies
\[
X\hbeta(y') = P_{X(\nul(D_{-\cA}))}(y') + c
\]
for all $y'$ in a neighborhood of $y$, where $c$ is a constant (not
depending on $y$). This again reveals an interesting connection to
linear regression, since it says that locally the generalized
lasso fit is a translation of the linear regression fit on $X$,
with the coefficients $\beta$ subject to $D_{-\cA}\beta=0$.
The same statements hold with the active set $\cA$ replaced by the
boundary set $\cB$ of an optimal subgradient.

We note that our results provide practically useful
estimates of degrees of freedom. For the lasso problem, we can use
$\rank(X_\cA)$ as an unbiased estimate of degrees of freedom, with
$\cA$ being the active set of a lasso solution. To emphasize what
has already been said, here we can actually choose any active set
(i.e., any solution), because all active sets give rise to the same
$\rank(X_\cA)$, except for $y$ in a set of measure zero. This is
important, since different algorithms for the lasso can produce
different solutions with different active sets. For
the generalized lasso problem, an unbiased estimate for degrees of
freedom is given by
$\dim(X(\nul(D_{-\cA})))=\rank(XP_{\nul(D_{-\cA})})$, where $\cA
$ is
the active set of a generalized lasso solution. This estimate is
the same, regardless of the choice of active set (i.e., choice of
solution), for almost every $y$. Hence any algorithm can be used to
compute a solution.

\begin{appendix}

\section{Proofs and technical arguments}
\label{app}

\subsection{\texorpdfstring{Proof of Lemma~\protect\ref{lem:nonexp}}{Proof of Lemma 1}}
\label{app:nonexp}

The proof relies on the fact that the projection $P_C(x)$ of $x \in
\R^n$ onto a closed convex set $C \subseteq\R^n$ satisfies
%
\begin{equation}
\label{eq:pfact}
\langle x-P_C(x), P_C(x)-u \rangle\geq0
\qquad\mbox{for any } u \in C.
\end{equation}

First, we prove the statement for the projection map. Note that
\begin{eqnarray*}
&&\|P_C(x)-P_C(y)\|_\ltwo^2
\\
 &&\qquad= \langle P_C(x)-x + y-P_C(y) + x-y,
P_C(x)-P_C(y) \rangle
\\
 &&\qquad= \langle P_C(x)-x, P_C(x)-P_C(y) \rangle+
\langle y-P_C(y), P_C(x)-P_C(y) \rangle
\\
&&\hphantom{\qquad=}{} + \langle x-y, P_C(x)-P_C(y) \rangle
\\
 &&\qquad\leq\langle x-y, P_C(x)-P_C(y) \rangle
\\
 &&\qquad\leq\|x-y\|_\ltwo\|P_C(x)-P_C(y)\|_\ltwo,
\end{eqnarray*}
where the first inequality follows from \eqref{eq:pfact}, and the
second is by Cauchy--Schwarz. Dividing both sides by
$\|P_C(x)-P_C(y)\|_\ltwo$ gives the result.

Now, for the residual map, the steps are similar.
\begin{eqnarray*}
&&\|(I-P_C)(x)-(I-P_C)(y)\|_\ltwo^2
\\
 &&\qquad=  \langle P_C(y)-P_C(x) + x-y,
x-P_C(x)+P_C(y)-y \rangle
\\
 &&\qquad= \langle P_C(y)-P_C(x), x-P_C(x) \rangle
+ \langle P_C(y)-P_C(x), P_C(y)-y \rangle
\\
&&\hphantom{\qquad=}{} + \langle x-y, x-P_C(x)+P_C(y)-y \rangle
\\
 &&\qquad\leq\langle x-y, x-P_C(x)+P_C(y)-y \rangle
\\
 &&\qquad\leq\|x-y\|_\ltwo\|(I-P_C)(x)-(I-P_C)(y)\|_\ltwo.
\end{eqnarray*}
Again the two inequalities are from \eqref{eq:pfact} and
Cauchy--Schwarz, respectively, and dividing both sides by
$\|(I-P_C)(x)-(I-P_C)(y)\|_\ltwo$ gives the result.

We have shown that $P_C$ and $I-P_C$ are Lipschitz (with constant
$1$); they are therefore continuous, and almost differentiability
follows from the standard proof of the fact that a Lipschitz function
is differentiable almost everywhere.

\subsection{\texorpdfstring{Proof of Lemma \protect\ref{lem:locaff}}{Proof of Lemma 2}}
\label{app:locaff}

We write $\cF$ to denote the set of faces of $C$. To each face $F \in
\cF$, there is an associated normal cone $N(F)$, defined as
\[
N(F) = \Bigl\{ x \in\R^n\dvtx F = \mathop\argmax_{y \in C} x\T y
\Bigr\}.
\]
The normal cone of $F$ satisfies $N(F)=P_C^{-1}(u)-u$
for any $u \in\relint(F)$. [We use $\relint(A)$ to
denote the relative interior of a set $A$, and $\relbd(A)$
to denote its relative boundary.]

Define the set
\[
\cS= \bigcup_{F \in\cF}
\bigl( \relint(F) + \relint(N(F)) \bigr).
\]
Because $C$ is a polyhedron, we have that $\dim(F)+\dim(N(F))=n$ for
each $F \in\cF$, and therefore each $U_F = \relint(F)+\relint(N(F))$
is an open set in~$\R^n$.

Now let $x \in\cS$.
We have $x \in U_F$ for some $F \in\cF$, and by construction
$P_C(U_F) = \relint(F)$. Furthermore, we claim that
projecting $x \in U_F$ onto $C$ is the same as projecting $x$ onto the
affine hull of $F$, that is, $P_C(U_F)=P_{\aff(F)}(U_F)$. Otherwise
there is some $y \in U_F$ with $P_C(y) \not= P_{\aff(F)}(y)$, and as
$\aff(F) \supseteq F$, this means that $\|y-P_{\aff(F)}(y)\|_\ltwo<
\|y-P_C(y)\|_\ltwo$. By definition of $\relint(F)$, there is some
$\alpha\in(0,1)$ such that $u=\alpha P_C(y) + (1-\alpha)P_{\aff(F)}
\in F$. But $\|y-u\|_\ltwo< \alpha\|y-P_C(y)\|_\ltwo+
(1-\alpha)\|y-P_{\aff(F)}(y)\|_\ltwo< \|y-P_C(y)\|_\ltwo$, which is a
contradiction. This proves the claim, and writing $\aff(F)=L+a$, we
have
\[
P_C(y) = P_L(y-a)+a \qquad\mbox{for } y\in U_F,
\]
as desired.\vadjust{\goodbreak}

It remains to show that $\cS^c=\R^n\setminus\cS$ has measure
zero. Note that $\cS^c$ contains points of the form $u+x$, where
either:
\begin{longlist}[(2)]
\item[(1)] $u \in\relbd(F), x \in N(F)$ for some $F$ with
$\dim(F)\geq1$; or
\item[(2)] $u \in\relint(F), x \in\relbd(N(F))$ for some $F\not=C$.
\end{longlist}
In the first type of points above, vertices are excluded because
$\relbd(F)=\varnothing$ when $F$ is a vertex. In the second type, $C$ is
excluded because $\relbd(N(C))=\varnothing$.
The lattice structure of $\cF$ tells us that for any face $F \in\cF$,
we can write $\relbd(F)=\bigcup_{G \in\cF, G \subsetneq F}
\relint(G)$. This, and the fact that the normal cones have the
opposite partial ordering as the faces, imply that points of the first
type above can be written as $u'+x'$ with $u' \in\relint(G)$ and $x'
\in N(G)$ for some $G \subsetneq F$. Note that actually we must have
$x' \in\relbd(N(G))$ because otherwise we would have $u'+x' \in
\cS$. Therefore it suffices to consider points of the second type
alone, and $\cS^c$ can be written as
\[
\cS^c = \bigcup_{F \in\cF, F \not= C}
\bigl(\relint(F)+\relbd(N(F)) \bigr).
\]
As $C$ is a polyhedron, the set $\cF$ of its faces is finite, and
$\dim(\relbd(N(F)))\leq n-\dim(F)-1$ for each $F \in\cF, F \not= C$.
Therefore $\cS^c$ is a finite union of sets of dimension $\leq n-1$,
and hence has measure zero.

\subsection{\texorpdfstring{Proof of Lemma \protect\ref{lem:lcact}}{Proof of Lemma 6}}
\label{app:lcact}

First some notation. For $S \subseteq\{1,\ldots, k\}$, define the
function $\pi_S \dvtx \R^k \rightarrow\R^{|S|}$ by $\pi_{S}(x) =
x_S$. So $\pi_S$ just extracts the coordinates in $S$.

Now let
\[
\cM= \bigcup_{\cE,s} \bigcup_{\cA\in Z(\cE)}
\bigl\{z \in\R^n \dvtx P_{[\pi_{-\cA}(\nul(X_\cE))]^\perp}
[(X_\cE)^+]_{(-\cA,\cdot)}
\bigl(z - (X_\cE\T)^+\lambda s \bigr) = 0 \bigr\}.
\]
The first union is taken over all possible subsets $\cE
\subseteq\{1,\ldots, p\}$ and all sign vectors $s\in
\{-1,1\}^{|\cE|}$; as for the second union, we define for a fixed
subset $\cE$
\[
Z(\cE) = \bigl\{ \cA\subseteq\cE\dvtx
P_{[\pi_{-\cA}(\nul(X_\cE))]^\perp}
[(X_\cE)^+]_{(-\cA,\cdot)} \not= 0 \bigr\}.
\]
Notice that $\cM$ is a finite union of affine subspace of dimension
$\leq n-1$, and hence has measure zero.

Let $y \notin\cM$, and let $\hbeta(y)$ be a lasso solution,
abbreviating $\cA=\cA(y)$ and $r=r(y)$ for the active set and active
signs. Also write $\cE=\cE(y)$ and $s=s(y)$ for the equicorrelation
set and equicorrelation signs of the fit. We know from~\eqref{eq:lassosol} that we can write
\[
\hbeta_{-\cE}(y) = 0\quad \mbox{and}\quad
\hbeta_\cE(y) = (X_\cE)^+ \bigl(y - (X_\cE\T)^+ \lambda s \bigr) + b,
\]
where $b \in\nul(X_\cE)$ is such that
\[
\hbeta_{\cE\setminus\cA}(y) =
[(X_\cE)^+]_{(-\cA,\cdot)}
\bigl(y - (X_\cE\T)^+ \lambda s \bigr) + b_{-\cA}
= 0.
\]
In other words,
\[
[(X_\cE)^+]_{(-\cA,\cdot)}
\bigl(y - (X_\cE\T)^+ \lambda s \bigr)
= -b_{-\cA}
\in \pi_{-\cA} (\nul(X_\cE)),
\]
so projecting onto the orthogonal complement of the linear
subspace\break $\pi_{-\cA} (\nul(X_\cE))$ gives zero,
\[
P_{[\pi_{-\cA}(\nul(X_\cE))]^\perp}
[(X_\cE)^+]_{(-\cA,\cdot)}
\bigl(y - (X_\cE\T)^+\lambda s \bigr) = 0.
\]
Since $y \notin\cM$, we know that
\[
P_{[\pi_{-\cA}(\nul(X_\cE))]^\perp}
[(X_\cE)^+]_{(-\cA,\cdot)} = 0,
\]
and finally, this can be rewritten as
%
\begin{equation}
\label{eq:colsp}
\col \bigl([(X_\cE)^+]_{(-\cA,\cdot)} \bigr)
\subseteq\pi_{-\cA}(\nul(X_\cE)).
\end{equation}

Consider defining, for a new point $y'$,
\[
\hbeta_{-\cE}(y') = 0 \quad\mbox{and}\quad
\hbeta_\cE(y') = (X_\cE)^+ \bigl(y' - (X_\cE\T)^+ \lambda s \bigr) + b',
\]
where $b' \in\nul(X_\cE)$, and is yet to be determined. Exactly as in
the proof of Lemma~\ref{lem:lcequi}, we know that
$X_\cE\T(y'-X\hbeta(y'))=\lambda s$, and
\mbox{$\|X_{-\cE}\T(y'-X\hbeta(y'))\|_\linf< \lambda$} for all $y' \in U_1$,
a neighborhood of $y$.

Now we want to choose $b'$ so that $\hbeta(y')$ has the
correct active set and active signs. For simplicity of notation,
first define the function $f \dvtx \R^n \rightarrow\R^{|\cE|}$,
\[
f(x) = (X_\cE)^+ \bigl(x - (X_\cE\T)^+ \lambda s \bigr).
\]
Equation \eqref{eq:colsp} implies that there is a
$b' \in\nul(X_\cE)$ such that
$b'_{-\cA}=-f_{-\cA}(y')$, hence
$\hbeta_{\cE\setminus\cA}(y')=0$. However,
we must choose $b'$ so that additionally
$\hbeta_i(y') \not= 0$ for $i \in\cA$ and
$\sign(\hbeta_\cA(y'))=r$. Write
\[
\hbeta_\cE(y')= \bigl(f(y')+b \bigr) + (b'-b).
\]
By the continuity of $f+b$, there exits a neighborhood of $U_2$ of
$y$ such that $f_i(y')+b_i\not=0$ for $i \in\cA$ and
$\sign(f_\cA(y')+b_\cA)=r$, for all $y' \in U_2$. Therefore
we only need to choose a vector $b' \in\nul(X_\cE)$, with
$b'_{-\cA}=-f_{-\cA}(y')$, such that $\|b'-b\|_\ltwo$ sufficiently
small. This can be achieved by applying the bounded inverse theorem,
which says that the bijective linear map $\pi_{-\cA}$ has a bounded
inverse (when considered a function from its row space to its column
space). Therefore there exists some $M>0$ such that for any $y'$,
there is a vector $b' \in\nul(X_\cE)$, $b'_{-\cA}=-f_{-\cA}(y')$,
with
\[
\|b'-b\|_\ltwo\leq M \|f_{-\cA}(y')-f_{-\cA}(y)\|_\ltwo.
\]
Finally, the continuity of $f_{-\cA}$ implies that
$\|f_{-\cA}(y')-f_{-\cA}(y)\|_\ltwo$ can be made sufficiently small by
restricting $y' \in U_3$, another neighborhood of $y$.

Letting $U=U_1\cap U_2\cap U_3$, we have shown that for any $y' \in
U$, there exists a lasso solution $\hbeta(y')$ with active set
$\cA(y')=\cA$ and active signs $r(y')=r$.

\subsection{\texorpdfstring{Proof of Lemma \protect\ref{lem:lcbound}}{Proof of Lemma 9}}
\label{app:lcbound}

Define the set
\begin{eqnarray*}
&&\cN= \bigcup_{\cB,s} \bigcup_{\cA\in Z(\cB)}
\bigl\{z \in\R^n\dvtx
P_{[D_{\cB\setminus\cA}(\nul(X) \cap\nul(D_{-\cB}))]^\perp}\cdot D_{\cB\setminus\cA} \bigl(X P_{\nul(D_{-\cB})}\bigr)^+
\\
&&\hspace*{165pt}{}\times
\bigl(z - \bigl(P_{\nul(D_{-\cB})} X\T\bigr)^+ D_\cB\T\lambda s \bigr) = 0
\bigr\}.
\end{eqnarray*}
The first union above is taken over all subsets
$\cB\subseteq\{1,\ldots, m\}$ and all sign vectors $s \in
\{-1,1\}^{|\cB|}$. The second union is taken over subsets $\cA
\subseteq Z(\cB)$, where
\[
Z(\cB) = \bigl\{ \cA\subseteq\cB\dvtx
P_{[D_{\cB\setminus\cA}(\nul(X) \cap\nul(D_{-\cB}))]^\perp}
D_{\cB\setminus\cA} \bigl(X P_{\nul(D_{-\cB})}\bigr)^+ \not= 0 \bigr\}.
\]
Since $\cN$ is a finite union of affine subspaces of dimension $\leq
n-1$, it has measure zero.

Now fix $y \notin\cN$, and let $(\gamma(y),\hbeta(y))$ be an
optimal pair, with boundary set $\cB=\cB(y)$, boundary
signs $s=s(y)$, active set $\cA=\cA(y)$, and active signs $r=r(y)$.
Starting from \eqref{eq:genlassokktb}, and plugging in for the fit in
terms of $\cB,s$, as in \eqref{eq:genlassofit} we can show that
\begin{eqnarray*}
\gamma_{-\cB}(y) &=& \lambda^{-1} (D_{-\cB}\T)^+
\bigl( X\T P_{\nul(P_{\nul(D_{-\cB})} X\T)} y
\\
&&\hphantom{\lambda^{-1} (D_{-\cB}\T)^+
\bigl(}{} +
\bigl(X\T\bigl(P_{\nul(D_{-\cB})} X\T\bigr)^+ - I \bigr) D_\cB\T\lambda s \bigr)
+ c,
\end{eqnarray*}
where $c \in\nul(D_{-\cB}\T)$. By \eqref{eq:genlassosol}, we
know that
\[
\hbeta(y) = \bigl(X P_{\nul(D_{-\cB})}\bigr)^+
\bigl(y - \bigl(P_{\nul(D_{-\cB})} X\T\bigr)^+ D_\cB\T\lambda s \bigr) + b,
\]
where $b \in\nul(X) \cap\nul(D_{-\cB})$. Furthermore,
\[
D_{\cB\setminus\cA}\hbeta(y) =
D_{\cB\setminus\cA} \bigl(X P_{\nul(D_{-\cB})}\bigr)^+
\bigl(y - \bigl(P_{\nul(D_{-\cB})} X\T\bigr)^+ D_\cB\T\lambda s \bigr)
+ D_{\cB\setminus\cA}b = 0,
\]
or equivalently,
\begin{eqnarray*}
&&D_{\cB\setminus\cA} \bigl(X P_{\nul(D_{-\cB})}\bigr)^+
\bigl(y - \bigl(P_{\nul(D_{-\cB})} X\T\bigr)^+ D_\cB\T\lambda s \bigr)
\\
&&\qquad =
-D_{\cB\setminus\cA}b
\in D_{\cB\setminus\cA} \bigl(\nul(X) \cap\nul(D_{-\cB})\bigr).
\end{eqnarray*}
Projecting onto the orthogonal complement of the
linear subspace\break $D_{\cB\setminus\cA}(\nul(X) \cap\nul(D_{-\cB}))$
therefore gives zero,
\[
P_{[D_{\cB\setminus\cA}(\nul(X) \cap\nul(D_{-\cB}))]^\perp}
D_{\cB\setminus\cA} \bigl(X P_{\nul(D_{-\cB})}\bigr)^+
\bigl(y - \bigl(P_{\nul(D_{-\cB})} X\T\bigr)^+ D_\cB\T\lambda s \bigr) = 0,
\]
and because $y \notin\cN$, we know that in fact
\[
P_{[D_{\cB\setminus\cA}(\nul(X) \cap\nul(D_{-\cB}))]^\perp}
D_{\cB\setminus\cA} \bigl(X P_{\nul(D_{-\cB})}\bigr)^+ = 0.
\]
This can be rewritten as
%
\begin{equation}
\label{eq:colsp2}
\col\bigl(D_{\cB\setminus\cA} \bigl(X P_{\nul(D_{-\cB})}\bigr)^+ \bigr)
\subseteq
D_{\cB\setminus\cA} \bigl(\nul(X) \cap\nul(D_{-\cB})\bigr).
\end{equation}

At a new point $y'$, consider defining $\gamma_\cB(y') = s$,
\begin{eqnarray*}
\gamma_{-\cB}(y')& =& \lambda^{-1} (D_{-\cB}\T)^+
\bigl( X\T P_{\nul(P_{\nul(D_{-\cB})} X\T)} y'
\\
&&\hphantom{\lambda^{-1} (D_{-\cB}\T)^+
\bigl(}{}+
\bigl(X\T\bigl(P_{\nul(D_{-\cB})} X\T\bigr)^+ - I \bigr) D_\cB\T\lambda s \bigr)
+ c,
\end{eqnarray*}
and
\[
\hbeta(y') = \bigl(X P_{\nul(D_{-\cB})}\bigr)^+
\bigl(y' - \bigl(P_{\nul(D_{-\cB})} X\T\bigr)^+ D_\cB\T\lambda s \bigr) + b',
\]
where $b' \in\nul(X)\cap\nul(D_{-\cB})$ is yet to be determined.
By construction, $\gamma(y')$ and $\hbeta(y')$ satisfy the
stationarity condition \eqref{eq:genlassokkt} at $y'$. Hence it
remains to show two parts: first, we must show that this pair
satisfies the subgradient condition \eqref{eq:genlassosg} at $y'$;
second, we must show this pair has boundary set \mbox{$\cB(y')=\cB$},
boundary signs $s(y')=s$, active set $\cA(y')=\cA$ and active signs
$r(y')=y$. Actually, it suffices to show the second part alone,
because the first part is then implied by the fact that
$\gamma(y)$ and $\hbeta(y)$ satisfy the subgradient condition at $y$.
Well, by the continuity of the function $f \dvtx \R^n \rightarrow
\R^{m-|\cB|}$,
\begin{eqnarray*}
f(x)& =& \lambda^{-1} (D_{-\cB}\T)^+
\bigl( X\T P_{\nul(P_{\nul(D_{-\cB})} X\T)} x
\\
&&\hphantom{\lambda^{-1} (D_{-\cB}\T)^+
\bigl(}{}+
\bigl(X\T\bigl(P_{\nul(D_{-\cB})} X\T\bigr)^+ - I \bigr) D_\cB\T\lambda s \bigr)
+ c,
\end{eqnarray*}
we have $\|\gamma_{-\cB}(y')\|_\linf< 1$ provided that $y' \in U_1$,
a neighborhood of $y$. This ensures that $\gamma(y')$ has boundary set
$\cB(y')=\cB$ and signs $s(y')=s$.

As for the active set and signs of $\hbeta(y')$, note first that
$D_{-\cB}\hbeta(y')=0$, following directly from the definition. Next,
define the function $g \dvtx \R^n \rightarrow\R^p$,
\[
g(x) = \bigl(X P_{\nul(D_{-\cB})}\bigr)^+
\bigl(x - \bigl(P_{\nul(D_{-\cB})} X\T\bigr)^+ D_\cB\T\lambda s \bigr),
\]
so $\hbeta(y') = g(y')+b'$. Equation \eqref{eq:colsp2} implies
that there is a vector $b' \in\nul(X)\cap\nul(D_{-\cB})$ such that
$D_{\cB\setminus\cA} b' = -D_{\cB\setminus\cA} g(y')$, which makes
$D_{\cB\setminus\cA}\hbeta(y')=0$. However, we still need to choose
$b'$ such that $D_i \hbeta(y') \not= 0$ for all $i \in\cA$ and
$\sign(D_\cA\hbeta(y'))=r$. To this end, write
\[
\hbeta(y') = \bigl(g(y')+b\bigr) + (b'-b).
\]
The continuity of $D_\cA g$ implies that there is a neighborhood
$U_2$ of $y$ such that $D_i g(y')+ D_i b \not=0$ for all $i \in\cA$
and $\sign(D_\cA g(y')+D_\cA b)=r$, for $y' \in U_2$. Since
\begin{eqnarray*}
|D_i\hbeta(y)|
&\geq&|D_ig(y')+D_ib| - |D_i(b'-b)|
\\
&\geq&|D_ig(y')+D_ib| - \|D\T\|_\ltwo\|b'-b\|_\ltwo,
\end{eqnarray*}
where $\|D\T\|_\ltwo$ is the operator norm of the $D\T$, we only need
to choose $b' \in\nul(X)\cap\nul(D_{-\cB})$ such that
$D_{\cB\setminus\cA} b' = -D_{\cB\setminus\cA} g(y')$, and such that
$\|b'-b\|_\ltwo$ is sufficiently small. This is possible by the
bounded inverse theorem applied to the linear map\vadjust{\goodbreak}
$D_{\cB\setminus\cA}$: when considered a function from its row space
to its column space, $D_{\cB\setminus\cA}$ is bijective and hence has
a bounded inverse. Therefore there is some $M>0$ such that for
any $y'$, there is a $b' \in\nul(X)\cap\nul(D_{-\cB})$ with
$D_{\cB\setminus\cA} b' = -D_{\cB\setminus\cA} g(y')$ and
\[
\|b'-b\|_\ltwo\leq M
\|D_{\cB\setminus\cA}g(y')-D_{\cB\setminus\cA}g(y)\|_\ltwo.
\]
The continuity of $D_{\cB\setminus\cA}g$ implies that the
right-hand side above can be made sufficiently small by
restricting $y' \in U_3$, a neighborhood of $y$.

With $U=U_1\cap U_2\cap U_3$, we have shown for that for $y' \in
U$, there is an optimal pair $(\gamma(y'),\hbeta(y'))$
with boundary set $\cB(y')=\cB$, boundary signs $s(y')=s$, active set
$\cA(y')=\cA$ and active signs $r(y')=r$.

\subsection{Dual problems}
\label{app:dual}

The dual of the lasso problem \eqref{eq:lasso} has appeared in many papers
in the literature; as far as we can tell, it was first considered
by \citet{homotopy}.
We start by rewriting problem \eqref{eq:lasso}~as
\[
\hbeta,\hat{z} \in\mathop\argmin_{\beta\in\R^p, z \in\R^n}
\half\|y-z\|_\ltwo^2 + \lambda\|\beta\|_\lone
\qquad\mbox{subject to } z=X\beta;
\]
then we write the Lagrangian
\[
\mathcal{L}(\beta,z,v) = \tfrac{1}{2}\|y-z\|_\ltwo^2 +
\lambda\|\beta\|_\lone+ v\T(z-X\beta),
\]
and we minimize $\mathcal{L}$ over $\beta,z$ to obtain the
dual problem
%
\begin{equation}
\label{eq:lassodual1}
\hv= \mathop\argmin_{v \in\R^n} \|y-v\|_\ltwo^2
\qquad\mbox{subject to } \|X\T v\|_\linf\leq\lambda.
\end{equation}
Taking the gradient of $\mathcal{L}$ with respect to
to $\beta,z$, and setting this equal to zero gives
%
\begin{eqnarray}
\label{eq:lassopd1}
\hv&=& y-X\hbeta,
\\
\label{eq:lassopd2}
X^T \hv&=& \lambda\gamma,
\end{eqnarray}
where $\gamma\in\R^p$ is a subgradient of the function $f(x)=\|x\|
_\lone$
evaluated at $x=\hbeta$. From \eqref{eq:lassodual1}, we can
immediately see
that the dual solution $\hv$ is the projection\vspace*{2pt} of $y$ onto the
polyhedron $C$
as in Lemma~\ref{lem:lassoproj},
and then \eqref{eq:lassopd1} shows that $X\hbeta=y-\hv$ is the
residual from
projecting $y$ onto $C$. Further, from~\eqref{eq:lassopd2},
we can define the equicorrelation set $\cE$ as
\[
\cE= \{ i \in\{1,\ldots, p\} \dvtx |X_i^T \hv| = \lambda \}.
\]
Noting that together \eqref{eq:lassopd1}, \eqref{eq:lassopd2} are exactly
the same as the KKT conditions~\eqref{eq:lassokkt}, \eqref
{eq:lassosg}, and all
of the arguments in Section~\ref{sec:lasso} involving the equicorrelation
set $\cE$ can be translated to this dual perspective.

There is a slightly different way to derive the lasso dual, resulting
in a~different (but of course, equivalent) formulation. We first rewrite
problem~\eqref{eq:lasso} as
\[
\hbeta,\hat{z} \in\mathop\argmin_{\beta\in\R^p, z \in\R^n}
\half\|y-X\beta\|_\ltwo^2 + \lambda\|z\|_\lone
\qquad\mbox{subject to } z=\beta,
\]
and by following similar steps to those above, we arrive at the dual problem
%
\begin{equation}
\label{eq:lassodual2}
\hv\in\mathop\argmin_{v \in\R^p} \bigl\|P_{\col(X)}y-(X^+)\T v\bigr\|_\ltwo^2
\qquad\mbox{subject to } \|v\|_\linf\leq\lambda, v \in\row(X).\hspace*{-35pt}
\end{equation}
Each dual solution $\hv$ (now no longer unique) satisfies
%
\begin{eqnarray}
\label{eq:lassopd3}
(X^+)\T\hv&=& P_{\col(X)}y - X\hbeta,
\\
\label{eq:lassopd4}
\hv&=& \lambda\gamma.
\end{eqnarray}
The dual problem \eqref{eq:lassodual2} and its relationship \eqref
{eq:lassopd3},
\eqref{eq:lassopd4} to the primal problem offer yet another viewpoint to
understand some of the results in Section~\ref{sec:lasso}.

For the generalized lasso problem, one might imagine that there are
three different dual problems, corresponding to the three different ways
of introducing an auxiliary variable $z$ into the generalized lasso criterion:
\begin{eqnarray*}
\hbeta,\hat{z}
&\in&\mathop\argmin_{\beta\in\R^p, z \in\R^n}
\half\|y-z\|_\ltwo^2 + \lambda\|D\beta\|_\lone
\qquad\mbox{subject to } z=X\beta;
\\[2pt]
\hbeta,\hat{z}
&\in&\mathop\argmin_{\beta\in\R^p, z \in\R^p}
\half\|y-X\beta\|_\ltwo^2 + \lambda\|Dz\|_\lone
\qquad\mbox{subject to } z=\beta;
\\[2pt]
\hbeta,\hat{z}
&\in&\mathop\argmin_{\beta\in\R^p, z \in\R^m}
\half\|y-X\beta\|_\ltwo^2 + \lambda\|z\|_\lone
\qquad\mbox{subject to } z=D\beta.
\end{eqnarray*}
However, the first two approaches above lead to Lagrangian functions that
cannot be minimized analytically over $\beta,z$. Only the third
approach yields
a dual problem in closed-form, as given by \citet{genlasso},
%
\begin{eqnarray}
\label{eq:genlassodual}
\hv\in\mathop\argmin_{v \in\R^m} \bigl\|P_{\col(X)}y-(X^+)\T D\T v\bigr\|_\ltwo^2\nonumber
\\[-7pt]
\\[-7pt]
\eqntext{\mbox{subject to } \|v\|_\linf\leq\lambda, D\T v \in\row(X).}
\end{eqnarray}
The relationship between primal and dual solutions is
%
\begin{eqnarray}
\label{eq:genlassopd1}
(X^+)\T D\T\hv&=& P_{\col(X)}y - X\hbeta,
\\[2pt]
\label{eq:genlassopd2}
\hv&= &\lambda\gamma,
\end{eqnarray}
where $\gamma\in\R^m$ is a subgradient of $f(x)=\|x\|_\lone$ evaluated
at $x=D\hbeta$. Directly from \eqref{eq:genlassodual} we can see that
$(X^+)\T D\T\hv$ is the projection of the point $y'=P_{\col(X)}y$ onto
the polyhedron
\[
K = \{ (X^+)\T D\T v \dvtx \|v\|_\linf\leq\lambda,
D\T v \in\row(X) \}.
\]
By \eqref{eq:genlassopd1}, the primal fit is $X\hbeta=(I-P_K)(y')$,
which can be rewritten as $X\hbeta=(I-P_C)(y')$ where
$C$ is the polyhedron from Lemma~\ref{lem:genlassoproj}, and finally
$X\hbeta=(I-P_C)(y)$ because $I-P_C$ is zero on $\nul(X^T)$.\vadjust{\goodbreak}
By \eqref{eq:genlassopd2}, we can define the boundary set $\cB$
corresponding to a particular dual solution $\hv$ as
\[
\cB= \bigl\{ i \in\{1,\ldots, m\} \dvtx |\hv_i| = \lambda \bigr\}.
\]
(This explains its name, as $\cB$ gives the coordinates of $\hv$
that are on the boundary of the box $\{v \dvtx \|v\|_\linf\leq\lambda\}$.)
As \eqref{eq:genlassopd1}, \eqref{eq:genlassopd2} are equivalent to
the KKT conditions \eqref{eq:genlassokkt}, \eqref{eq:genlassosg}
[following from rewriting \eqref{eq:genlassopd1} using $D\T\hv\in
\row(X)$],
the results in Section~\ref{sec:genlasso} on the boundary set $\cB$
can all be derived from this dual setting.
\end{appendix}



\printaddresses


\begin{thebibliography}{21}

\bibitem[\protect\citeauthoryear{Chen, Donoho and Saunders}{1998}]{bp}
\begin{barticle}[mr]
\bauthor{\bsnm{Chen},~\bfnm{Scott~Shaobing}\binits{S.~S.}},
  \bauthor{\bsnm{Donoho},~\bfnm{David~L.}\binits{D.~L.}} \AND
  \bauthor{\bsnm{Saunders},~\bfnm{Michael~A.}\binits{M.~A.}}
(\byear{1998}).
\btitle{Atomic decomposition by basis pursuit}.
\bjournal{SIAM J. Sci. Comput.}
\bvolume{20}
\bpages{33--61}.
\bid{doi={10.1137/S1064827596304010}, issn={1064-8275}, mr={1639094}}
\bptok{imsref}%
\end{barticle}
\endbibitem

\bibitem[\protect\citeauthoryear{Dossal et~al.}{2011}]{dfl1}
\begin{bmisc}[auto:STB|2012/04/30|08:06:40]
\bauthor{\bsnm{Dossal},~\bfnm{C.}\binits{C.}},
  \bauthor{\bsnm{Kachour},~\bfnm{M.}\binits{M.}},
  \bauthor{\bsnm{Fadili},~\bfnm{J.}\binits{J.}},
  \bauthor{\bsnm{Peyre},~\bfnm{G.}\binits{G.}} \AND
  \bauthor{\bsnm{Chesneau},~\bfnm{C.}\binits{C.}}
(\byear{2011}).
\bhowpublished{The degrees of freedom of the lasso for general design matrix.
  Available at arXiv:\arxivurl{1111.1162}}.
\bptok{imsref}%
\end{bmisc}
\endbibitem

\bibitem[\protect\citeauthoryear{Efron}{1986}]{bradbiased}
\begin{barticle}[mr]
\bauthor{\bsnm{Efron},~\bfnm{Bradley}\binits{B.}}
(\byear{1986}).
\btitle{How biased is the apparent error rate of a prediction rule?}
\bjournal{J. Amer. Statist. Assoc.}
\bvolume{81}
\bpages{461--470}.
\bid{issn={0162-1459}, mr={0845884}}
\bptok{imsref}%
\end{barticle}
\endbibitem

\bibitem[\protect\citeauthoryear{Efron et~al.}{2004}]{lars}
\begin{barticle}[mr]
\bauthor{\bsnm{Efron},~\bfnm{Bradley}\binits{B.}},
  \bauthor{\bsnm{Hastie},~\bfnm{Trevor}\binits{T.}},
  \bauthor{\bsnm{Johnstone},~\bfnm{Iain}\binits{I.}} \AND
  \bauthor{\bsnm{Tibshirani},~\bfnm{Robert}\binits{R.}}
(\byear{2004}).
\btitle{Least angle regression (with discussion, and a rejoinder by the
  authors)}.
\bjournal{Ann. Statist.}
\bvolume{32}
\bpages{407--499}.\break
\bid{doi={10.1214/009053604000000067}, issn={0090-5364}, mr={2060166}}
\bptok{imsref}%
\end{barticle}
\endbibitem

\bibitem[\protect\citeauthoryear{Fan and Li}{2001}]{scad}
\begin{barticle}[mr]
\bauthor{\bsnm{Fan},~\bfnm{Jianqing}\binits{J.}} \AND
  \bauthor{\bsnm{Li},~\bfnm{Runze}\binits{R.}}
(\byear{2001}).
\btitle{Variable selection via nonconcave penalized likelihood and its oracle
  properties}.
\bjournal{J. Amer. Statist. Assoc.}
\bvolume{96}
\bpages{1348--1360}.
\bid{doi={10.1198/016214501753382273}, issn={0162-1459}, mr={1946581}}
\bptok{imsref}%
\end{barticle}
\endbibitem

\bibitem[\protect\citeauthoryear{Gr{\"u}nbaum}{2003}]{grunbaum}
\begin{bbook}[mr]
\bauthor{\bsnm{Gr{\"u}nbaum},~\bfnm{Branko}\binits{B.}}
(\byear{2003}).
\btitle{Convex Polytopes},
\bedition{2nd} ed.
\bseries{Graduate Texts in Mathematics}
\bvolume{221}.
\bpublisher{Springer}, \baddress{New York}.
\bid{doi={10.1007/978-1-4613-0019-9}, mr={1976856}}
\bptok{imsref}%
\end{bbook}
\endbibitem

\bibitem[\protect\citeauthoryear{Hastie and Tibshirani}{1990}]{gam}
\begin{bbook}[mr]
\bauthor{\bsnm{Hastie},~\bfnm{T.~J.}\binits{T.~J.}} \AND
  \bauthor{\bsnm{Tibshirani},~\bfnm{R.~J.}\binits{R.~J.}}
(\byear{1990}).
\btitle{Generalized Additive Models}.
\bseries{Monographs on Statistics and Applied Probability}
\bvolume{43}.
\bpublisher{Chapman \& Hall}, \baddress{London}.
\bid{mr={1082147}}
\bptok{imsref}%
\end{bbook}
\endbibitem

\bibitem[\protect\citeauthoryear{Loubes and Massart}{2004}]{loubes}
\begin{bmisc}[auto:STB|2012/04/30|08:06:40]
\bauthor{\bsnm{Loubes},~\bfnm{J.~M.}\binits{J.~M.}} \AND
  \bauthor{\bsnm{Massart},~\bfnm{P.}\binits{P.}}
(\byear{2004}).
\bhowpublished{Dicussion to ``{L}east angle regression.'' \textit{Ann.
  Statist.} \textbf{32} 460--465}.
\bptok{imsref}%
\end{bmisc}
\endbibitem

\bibitem[\protect\citeauthoryear{Mallows}{1973}]{mallows}
\begin{barticle}[auto:STB|2012/04/30|08:06:40]
\bauthor{\bsnm{Mallows},~\bfnm{C.}\binits{C.}}
(\byear{1973}).
\btitle{Some comments on {$C_p$}}.
\bjournal{Technometrics}
\bvolume{15}
\bpages{661--675}.
\bptok{imsref}%
\end{barticle}
\endbibitem

\bibitem[\protect\citeauthoryear{Meyer and Woodroofe}{2000}]{meyerwood}
\begin{barticle}[mr]
\bauthor{\bsnm{Meyer},~\bfnm{Mary}\binits{M.}} \AND
  \bauthor{\bsnm{Woodroofe},~\bfnm{Michael}\binits{M.}}
(\byear{2000}).
\btitle{On the degrees of freedom in shape-restricted regression}.
\bjournal{Ann. Statist.}
\bvolume{28}
\bpages{1083--1104}.
\bid{doi={10.1214/aos/1015956708}, issn={0090-5364}, mr={1810920}}
\bptok{imsref}%
\end{barticle}
\endbibitem

\bibitem[\protect\citeauthoryear{Osborne, Presnell and
  Turlach}{2000}]{homotopy}
\begin{barticle}[mr]
\bauthor{\bsnm{Osborne},~\bfnm{Michael~R.}\binits{M.~R.}},
  \bauthor{\bsnm{Presnell},~\bfnm{Brett}\binits{B.}} \AND
  \bauthor{\bsnm{Turlach},~\bfnm{Berwin~A.}\binits{B.~A.}}
(\byear{2000}).
\btitle{On the {LASSO} and its dual}.
\bjournal{J.~Comput. Graph. Statist.}
\bvolume{9}
\bpages{319--337}.
\bid{doi={10.2307/1390657}, issn={1061-8600}, mr={1822089}}
\bptok{imsref}%
\end{barticle}
\endbibitem

\bibitem[\protect\citeauthoryear{Rosset, Zhu and Hastie}{2004}]{boostpath}
\begin{barticle}[mr]
\bauthor{\bsnm{Rosset},~\bfnm{Saharon}\binits{S.}},
  \bauthor{\bsnm{Zhu},~\bfnm{Ji}\binits{J.}} \AND
  \bauthor{\bsnm{Hastie},~\bfnm{Trevor}\binits{T.}}
(\byear{2004}).
\btitle{Boosting as a regularized path to a maximum margin classifier}.
\bjournal{J. Mach. Learn. Res.}
\bvolume{5}
\bpages{941--973}.
\bid{issn={1532-4435}, mr={2248005}}
\bptnote{check year}%
\bptok{imsref}%
\end{barticle}
\endbibitem

\bibitem[\protect\citeauthoryear{Schneider}{1993}]{schneider}
\begin{bbook}[mr]
\bauthor{\bsnm{Schneider},~\bfnm{Rolf}\binits{R.}}
(\byear{1993}).
\btitle{Convex Bodies: The {B}runn--{M}inkowski Theory}.
\bseries{Encyclopedia of Mathematics and Its Applications}
\bvolume{44}.
\bpublisher{Cambridge Univ. Press}, \baddress{Cambridge}.
\bid{doi={10.1017/CBO9780511526282}, mr={1216521}}
\bptok{imsref}%
\end{bbook}
\endbibitem

\bibitem[\protect\citeauthoryear{Stein}{1981}]{stein}
\begin{barticle}[mr]
\bauthor{\bsnm{Stein},~\bfnm{Charles~M.}\binits{C.~M.}}
(\byear{1981}).
\btitle{Estimation of the mean of a multivariate normal distribution}.
\bjournal{Ann. Statist.}
\bvolume{9}
\bpages{1135--1151}.
\bid{issn={0090-5364}, mr={0630098}}
\bptok{imsref}%
\end{barticle}
\endbibitem

\bibitem[\protect\citeauthoryear{Tibshirani}{1996}]{lasso}
\begin{barticle}[mr]
\bauthor{\bsnm{Tibshirani},~\bfnm{Robert}\binits{R.}}
(\byear{1996}).
\btitle{Regression shrinkage and selection via the lasso}.
\bjournal{J. Roy. Statist. Soc. Ser. B}
\bvolume{58}
\bpages{267--288}.
\bid{issn={0035-9246}, mr={1379242}}
\bptok{imsref}%
\end{barticle}
\endbibitem

\bibitem[\protect\citeauthoryear{Tibshirani}{2011}]{ryanphd}
\begin{bmisc}[auto:STB|2012/04/30|08:06:40]
\bauthor{\bsnm{Tibshirani},~\bfnm{R.~J.}\binits{R.~J.}}
(\byear{2011}).
\bhowpublished{The solution path of the generalized lasso, Ph.D. thesis, Dept.
  Statistics, Stanford Univ.}
\bptok{imsref}%
\end{bmisc}
\endbibitem

\bibitem[\protect\citeauthoryear{Tibshirani and Taylor}{2011}]{genlasso}
\begin{barticle}[mr]
\bauthor{\bsnm{Tibshirani},~\bfnm{Ryan~J.}\binits{R.~J.}} \AND
  \bauthor{\bsnm{Taylor},~\bfnm{Jonathan}\binits{J.}}
(\byear{2011}).
\btitle{The solution path of the generalized lasso}.
\bjournal{Ann. Statist.}
\bvolume{39}
\bpages{1335--1371}.
\bid{doi={10.1214/11-AOS878}, issn={0090-5364}, mr={2850205}}
\bptok{imsref}%
\end{barticle}
\endbibitem

\bibitem[\protect\citeauthoryear{Vaiter et~al.}{2011}]{analreg}
\begin{bmisc}[auto:STB|2012/04/30|08:06:40]
\bauthor{\bsnm{Vaiter},~\bfnm{S.}\binits{S.}},
  \bauthor{\bsnm{Peyre},~\bfnm{G.}\binits{G.}},
  \bauthor{\bsnm{Dossal},~\bfnm{C.}\binits{C.}} \AND
  \bauthor{\bsnm{Fadili},~\bfnm{J.}\binits{J.}}
(\byear{2011}).
\bhowpublished{Robust sparse analysis regularization. Available at
  arXiv:\arxivurl{1109.6222}}.
\bptok{imsref}%
\end{bmisc}
\endbibitem

\bibitem[\protect\citeauthoryear{Zou and Hastie}{2005}]{enet}
\begin{barticle}[mr]
\bauthor{\bsnm{Zou},~\bfnm{Hui}\binits{H.}} \AND
  \bauthor{\bsnm{Hastie},~\bfnm{Trevor}\binits{T.}}
(\byear{2005}).
\btitle{Regularization and variable selection via the elastic net}.
\bjournal{J. R. Stat. Soc. Ser. B Stat. Methodol.}
\bvolume{67}
\bpages{301--320}.
\bid{doi={10.1111/j.1467-9868.2005.00503.x}, issn={1369-7412}, mr={2137327}}
\bptok{imsref}%
\end{barticle}
\endbibitem

\bibitem[\protect\citeauthoryear{Zou, Hastie and Tibshirani}{2007}]{lassodf}
\begin{barticle}[mr]
\bauthor{\bsnm{Zou},~\bfnm{Hui}\binits{H.}},
  \bauthor{\bsnm{Hastie},~\bfnm{Trevor}\binits{T.}} \AND
  \bauthor{\bsnm{Tibshirani},~\bfnm{Robert}\binits{R.}}
(\byear{2007}).
\btitle{On the ``degrees of freedom'' of the lasso}.
\bjournal{Ann. Statist.}
\bvolume{35}
\bpages{2173--2192}.
\bid{doi={10.1214/009053607000000127}, issn={0090-5364}, mr={2363967}}
\bptok{imsref}%
\end{barticle}
\endbibitem

\end{thebibliography}
\end{document}